\documentclass[11pt]{article}
\usepackage{psboxit,amssymb,amstex}

\PScommands
\newcommand{\bA}{{\bf{A}}}
\newcommand{\bB}{{\bf{B}}}
\newcommand{\bb}{{\bf{b}}}

\newcommand{\be}{{\bf{e}}}
\newcommand{\bH}{{\bf{H}}}

\newcommand{\bG}{{\bf{G}}}
\newcommand{\bI}{{\bf{I}}}
\newcommand{\bJ}{{\bf{J}}}

\newcommand{\bk}{{\bf{k}}}

\newcommand{\bR}{{\bf{R}}}

\newcommand{\bS}{{\bf{S}}}

\newcommand{\bT}{{\bf{T}}}

\newcommand{\bu}{{\bf{u}}}
\newcommand{\bU}{{\bf{U}}}
\newcommand{\bV}{{\bf{V}}}

\newcommand{\bx}{{\bf{x}}}

\newcommand{\by}{{\bf{y}}}
\newcommand{\bZ}{{\bf{Z}}}
\newcommand{\bz}{{\bf{z}}}
\newcommand{\bv}{{\bf{v}}}
\newcommand{\bW}{{\bf{W}}}
\newcommand{\bw}{{\bf{w}}}
\newcommand{\nrm}{{\rule[-1mm]{.6mm}{4mm}}}

\newcommand{\scalH}{{ \mbox{\raisebox{-.4ex}{\scriptsize ${\cal H}$}} }}
\newcommand{\scalV}{{ \mbox{\raisebox{-.4ex}{\scriptsize ${\cal V}$}} }}
\newcommand{\scalU}{{ \mbox{\raisebox{-.4ex}{\scriptsize ${\cal U}$}} }}
\numberwithin{equation}{section}
\newtheorem{theorem}{Theorem}[section] 
\newtheorem{corollary}[theorem]{Corollary}
\newtheorem{lemma}[theorem]{Lemma}
  \title{Galerkin Eigenvector Approximations\thanks{This work was 
  supported under the auspices of AFOSR Grant F49620-96-1-0329} }
  \author{Christopher Beattie\thanks{Department of Mathematics,
   Virginia Polytechnic Institute and State University, 
   Blacksburg, VA  24061  USA.}  } 

\begin{document}
\maketitle
\parbox[t]{4in}{How close are Galerkin eigenvectors to the best 
 approximation available out of the trial subspace  ? 
Under a variety of conditions the Galerkin method gives an approximate 
eigenvector that approaches asymptotically the projection of the exact eigenvector
onto the trial subspace -- and this occurs more rapidly than the 
underlying rate of convergence of the approximate eigenvectors.   
Both orthogonal-Galerkin and 
Petrov-Galerkin methods are considered here with a special emphasis 
on nonselfadjoint problems.   Consequences for the numerical treatment of
elliptic PDEs discretized either with finite element methods or with  spectral methods
are discussed and an application to Krylov subspace methods
 for large scale matrix eigenvalue problems is presented. 
New lower bounds to the $sep$ of a pair of operators are developed as 
well.
}

\vspace{4mm}

\vspace{4mm}

\section{Introduction}
  Consider the eigenvalue problem for  a linear operator $A$:
 \begin{equation}
 \begin{array}{c}
 \mbox{Find $\lambda \in {\Bbb C}$ and $ \hat{v}\neq 0$ so that}\\
 	A\hat{v}=\lambda \hat{v} 	
 	\label{op_eig}
 \end{array}
 \end{equation}
 We seek a family of approximations  
 $\{ \lambda_{h},\ \hat{v}_{h}\}_{h>0}$ using the Galerkin method.

 The Galerkin method approximates the operator $A$ with a finite rank 
 (or in finite dimensions, a much lower rank) operator, $A_h$ --  the ``projection'' of $A$, 
 that samples the action of $A$ on a subspace. The solution to 
 (\ref{op_eig}) is then approximated with a matrix eigenvalue problem 
 associated with $A_h$.
  
This work focuses on one 
particular bit of Galerkin folklore -- ``the Galerkin method yields 
an approximate eigenvector for $A$ that is essentially the projection of the 
exact eigenvector $\hat{v}$ onto the trial subspace.''  We discover 
 that this statement is correct under some mild 
conditions if\, 1)  ``essentially'' is taken to mean ``asymptotically,'' and 
\, 2) the projection involved is intrinsic to the 
Galerkin method and may be either orthogonal or oblique depending on 
how the discretization is organized and what point of view is taken.  Although more generality is possible, we restrict 
ourselves to a Hilbert space setting -- specific assumptions are
found in Section 2.

\setlength{\unitlength}{1.3mm}
\begin{picture}(50,60)(0,-20)
\thinlines
\put(0,0){\line(4,3){32}}
\put(32,24){\line(1,0){50}}
\put(17,5){\circle*{1}}
\put(2,0){\line(3,1){72}}
\put(17,4){\makebox(0,0)[tl]{$0$}}
\put(36,33){\makebox(0,0){$\hat{v}$}}
\multiput(35,32)(.7,-4.2){5}{\line(1,-6){.5}}
\multiput(35,32)(0,-2){9}{\line(0,-1){1}}
\put(31,15){\makebox(0,0){$P_{h}\hat{v}$}}
\put(38,9){\makebox(0,0){$\hat{v}_{h}$}}
\put(36.5,13.5){\makebox(0,0){?}}
\thicklines
\put(17,5){\vector(2,1){18}}
\put(17,5){\vector(3,1){21.316}}
\put(17,5){\vector(2,3){18}}
\put(0,0){\line(1,0){50}}
\put(49,-1){\makebox(0,0)[tl]{${\cal S}_{1,h}$}}
\put(82,24){\line(-4,-3){32}}
\put(68,21){\makebox(0,0)[tl]{${\cal U}_{h}$}}
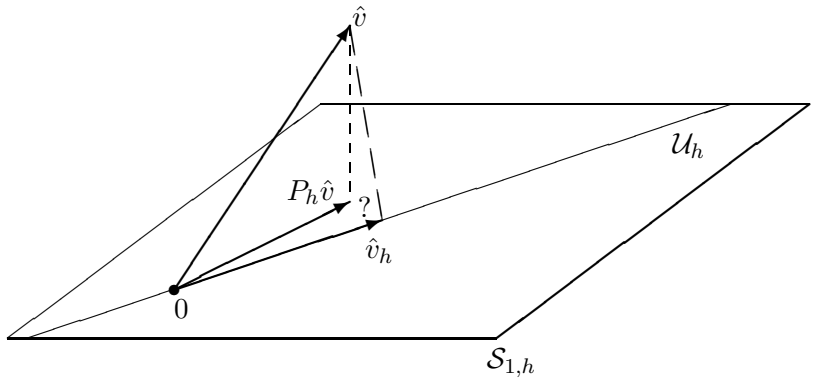
\put(10,-11){\shortstack{Figure 1.  How close is the
approximate eigenvector $\hat{v}_{h}$  \\ 
to the projected exact eigenvector $P_{h}\hat{v}$ ?}}
\end{picture}  

The basic features of Galerkin methods that play a role in 
our analysis are reviewed in Section 3. Of particular note here is 
that discussion is not restricted to self-adjoint problems.
Section 4 provides analysis for the simplest case -- when 
$A$ is a bounded operator.  The case where $A$ is an
unbounded operator is considered from two vantage points in the next 
two sections: with respect to the 
``energy'' norm in Section 5 where a discussion of consequences for the 
finite element method on elliptic problems may be found; and with respect to the underlying Hilbert space 
norm in Section 6 where an elliptic problem dicretized using a spectral method 
is discussed.  The final section presents an application to the 
biorthogonal Lanczos method and Arnoldi method.

\section{Setting of the Problem}   \label{setting}
\subsection{Operators defined via quadratic forms}
Although eigenvalue problems are most  naturally  posed 
for linear operators, the operators themselves are often difficult to 
specify fully -- particularly with regard to the operator's precise 
domain of definition.  It is often easier to characterize an operator in 
terms of a quadratic form that is naturally associated with it.  This 
approach usually leads spontaneously to the appropriate choice of 
underlying Hilbert spaces.  The reader may refer to the excellent 
tract of Kato \cite{Kato} for background material on quadratic forms.

Let ${\cal H}$ be a complex separable Hilbert space with inner 
product\footnote{Inner products and sesquilinear forms are 
conjugate linear in the first argument.} and norm denoted by  
$\langle\ \cdot\ ,\ \cdot\ \rangle_\scalH$ and $\|\cdot\  
\|_\scalH$, respectively.  Let $a(\cdot ,\ \cdot )$  be a
closed sectorial sesquilinear form, densely defined  in
${\cal H}$.  ``Sectorial'' means that 
\begin{equation}
	\begin{array}{c}
			\Re e\  a(v,v)  \geq  \alpha \|v\|^{2}_\scalH \\
  |\Im m\  a(v,v)|  \leq  M (\Re  e\   a(v,v) - \alpha \|v\|^{2}_\scalH)
	\end{array}
\label{semibnd}
\end{equation}
for all $ v \in {Dom(a)}$ and some fixed $\alpha,\ M >0$.
Following the notation of Kato \cite{Kato}, define symmetric sesquilinear forms associated 
with $a$,
\begin{equation}
	\begin{array}{c}
		[\Re e\  a](w,v)  =  \frac{1}{2} (a(w,v)+\overline{a(v,w)}) \\
		{[}\Im m\  a{]}(w,v)  =  \frac{1}{2\imath} (a(w,v)-\overline{a(v,w)}),
	\end{array}
	\label{ReIm}
\end{equation}
so that $a(w,v)= [\Re e\  a](w,v) +  \imath\ [\Im m\  a](w,v)$.  Notice that 
(\ref{semibnd}) implies
\begin{equation*}
	\alpha \|v\|^{2}_\scalH \leq [\Re e\  a](v,v) \leq |a(v,v)| 
	\leq \sqrt{1+M^{2}}\, [\Re e\  a](v,v),
\end{equation*}
so $[\Re e\  a] $ is a closed, symmetric, positive-definite, sesquilinear form
that induces an inner product on $Dom(a)$ with respect to which $Dom(a)$ is
a Hilbert space.   $a(u,v)$ is then a bounded sesquilinear form on this 
Hilbert space.  Furthermore, there is a  closed 
operator, $C_{a}$, densely defined on ${\cal H}$ so that 
$Dom(C_{a})=Dom(a)$ and $[\Re e\  a](v,v)=\|C_{a}v\|_\scalH^{2}$ 
(cf. \cite{Kato}, p. 331). 

Suppose now that ${\cal V}=Dom(a)$ is equipped with 
an inner product $\langle\ \cdot\ ,
\ \cdot\ \rangle_\scalV$ equivalent to $[\Re e\  a] $.   The Hilbert space
${\cal V}$ is continuously and densely imbedded in ${\cal H}$ and 
we may assume without loss of generality that 
$ \|u\|_\scalH \leq \|u\|_\scalV$  for 
 all $u \in  {\cal V}$. 

 Observe that any ${\cal H}$-bounded linear functional on  ${\cal H}$ 
 may be viewed immediately as   the extension of some
  ${\cal V}$-bounded linear functional on ${\cal V}$, so letting
  ${\cal V}'$ denote the dual space of ${\cal V}$, the
  imbedding  ${\cal V} \hookrightarrow {\cal H} $
may be extended to a Gelfand triple (see e.g., \cite{wloka}) ${\cal V} \hookrightarrow {\cal H} 
\hookrightarrow {\cal V}'$ with the norm on  ${\cal V}'$  defined by
$$
\|v\|_{\cal V'}=\sup_{w\in \cal V}\frac{|\langle w,v\rangle_{\cal 
H}|}{\|w\|_\scalV}.
$$
 The Cauchy-Schwartz inequality yields
  $ \| v  \|_{\cal V'} \leq\| v  \|_\scalH$
  for all $v \in  {\cal H}$.   

Under the hypotheses given,  Kato's first representation theorem 
(\cite{Kato}, p. 322) 
guarantees the existence of  a closed m-sectorial operator, 
$A$, defined on 
$$
Dom(A)=\left\{ v\in {\cal V}\left| \ |a(w,v)|\leq m_{v} \|w\|_{\cal 
H} \mbox{ for all $w\in Dom(a)$}\right. \right\}
$$
where $m_{v}$ is independent of $w$ but will generally depend on $v$).
Then 
$$
a(u,v)=\langle u,Av\rangle_\scalH
$$
for all $v \in Dom(A)$ and $ u \in {\cal V} $.  Furthermore, 
there is a closed operator $B_{a}$ with 
$Dom(B_{a})=Dom(C_{a})=Dom(a)$, so that $A$ may be decomposed 
as $A=B_{a}^{*}C_{a}$  (see e.g., \cite{Kato} p. 337). 
 ``$*$'' denotes the ${\cal H}$-adjoint. 

Since $Dom(A)$ 
is dense in ${\cal V}$ (with respect to the ${\cal V}$-norm) and
$$
	|a(u, v)|  \le  c \|u\|_\scalV \|v\|_\scalV,
$$
 we may calculate for any $v \in Dom(A)$
$$
\|Av\|_{\cal V'}=\sup_{u\in \cal V}\frac{|\langle u,Av\rangle_{\cal 
H}|}{\|u\|_\scalV} \leq c \|v\|_\scalV
$$
Thus $A$ may be extended by continuity to a bounded linear 
transformation from ${\cal V}$ to ${\cal V}'$.
Furthermore, 
$$
\|Av\|_{\cal V'}=\sup_{u\in \cal V}\frac{|a(u,v)|}{\|u\|_\scalV} \geq 
\frac{|a(v,v)|}{\|v\|_\scalV} \geq \frac{[\Re e\  a](v,v)}{\|v\|_\scalV} 
\geq {\hat \alpha} \|v\|_\scalV,
$$
thus $A$ is an isomorphism of ${\cal V}$ onto ${\cal V}'$ with a  
bounded inverse, $T=A^{-1}$ mapping ${\cal V}' $ back to $ {\cal V}$.
$T$ may be defined alternatively for each $v\in {\cal V}'$, so that $Tv\in 
{\cal V}$ solves 
\begin{equation}
	\langle u,\, v\rangle_\scalH=a(u,Tv)
	\label{BTdef}
\end{equation}
for all $u\in {\cal V}$.
 $A$ maps vectors in $Dom(A)$ to ${\cal H}$,  whereas $T$ maps vectors in ${\cal H}\subset 
{\cal V}'$ back to $Dom(A) \subset {\cal V}$. 

\subsection{The gap} 
Given two  closed subspaces, ${\cal  M}$ and ${\cal N}$ of ${\cal H}$,  the proximity 
of one to the other is measured in terms of the 
{\em containment gap } (or just {\em gap}\footnote{Kato defines the {\it gap} as 
$\max [\delta({\cal  M},{\cal  N}),\ \delta({\cal  N},{\cal  M})]$}), 
which we define   as
      $$
      \delta_\scalH({\cal  M},{\cal  N}) 
	 =  
	\sup_{ x\in{\cal  M}}\inf_{\,y\in{\cal  N}} 
	   \frac{\| y -  x \|_\scalH}{\| x\|_\scalH}= \sin(\Theta_{max}({\cal M},{\cal N})). 
	   $$ 
$\Theta_{max}({\cal M},{\cal  N})$ is the largest canonical angle between ${\cal  M}$ 
and a ``closest'' subspace $\widehat {\cal  N}$ of ${\cal  N}$ having dimension equal to 
$\dim{\cal M}$. Notice that if $\dim {\cal  N}< \dim {\cal  M}$ then $\delta_\scalH({\cal M},{\cal  N})=1$ and
$\delta_\scalH({\cal  M},{\cal  N})=0$ if and only if 
${\cal  M}\subset {\cal  N}$.   
 If $\dim {\cal  N} = \dim {\cal  M}$ 
then $\delta_\scalH({\cal  M},{\cal  N})=\delta_\scalH({\cal  N},{\cal  M})$.  
Conversely, if both $\delta_\scalH({\cal  M},{\cal  N})<1$ and $\delta_\scalH({\cal  N},{\cal  
M})<1$ then $\delta_\scalH({\cal  M},{\cal  N})=\delta_\scalH({\cal  N},{\cal  M})$ 
and $\dim {\cal  N} = \dim {\cal  M}$. 

The  gap can be  expressed directly as the norm of a 
composition of projections, so that if $\Pi_{{\cal  M}}$ and $\Pi_{{\cal  N}}$ denote 
${\cal  H}$-orthogonal projections onto ${\cal  M}$ and ${\cal N} $, 
respectively, then $\delta_\scalH({\cal  M},{\cal  N})=
\|(I-\Pi_{{\cal N}})\Pi_{{\cal  M}}\|_\scalH$. 

  If ${\cal  M}$ and ${\cal N} $ are closed subspaces of 
${\cal  V}$, we have the completely analogous definition of gap 
relative to ${\cal V}$:
 $$
      \delta_{\cal  V}({\cal  M},{\cal  N}) 
	 =  
	\sup_{ x\in{\cal  M}}\inf_{\,y\in{\cal  N}} 
	   \frac{\| y -  x \|_\scalV}{\| x\|_\scalV}. 
	   $$ 
	   
If $A$ and $B$ are closed operators in ${\cal  H}$ the gap between  $A$ and $B$
is defined as the gap between their graphs, 
considered as subspaces of ${\cal H}\times{\cal  H}$:
  $$
      \delta_\scalH(A, B) 
	 =  \sup_{x\in Dom(A)} \inf_{\,y\in Dom(B)} 
	   \frac{\| x -  y \|_\scalH+\| Ax -  By \|_\scalH }{ \|x\|_\scalH+\|Ax\|_\scalH}. 
	   $$

\subsection{The eigenvalue problem}
Our focus  rests on the (weakly posed) eigenvalue problem for  $a$:
 \begin{equation}
 \begin{array}{c}
 \mbox{Find $\lambda$ and $0\neq \hat{v}\in Dom(a)$ so that}\\
 	a(w,\hat{v})=\lambda \langle w,\hat{v} \rangle_\scalH 	
 	\mbox{  for all $w\in {Dom(a)}$.}
 	\label{var_eig}
 \end{array}
 \end{equation}
 Note that  $\{\lambda,\ \hat{v}\}$  
is an eigenpair for (\ref{var_eig}) if and 
only if $ \hat{v} \in Dom(A)$ and $\{\lambda,\ \hat{v}\}$ is an eigenpair 
for the operator $A$; or equivalently 
when $\lambda \neq 0$, if
 $\{\lambda^{-1},\ \hat{v}\}$ is an eigenpair for the operator $T$.

 Denote the {\it resolvent set} of $A$ by 
$$
\rho(A)=\left\{z\in {\Bbb C}\left| z-A \mbox{ has a bounded inverse on }
{\cal H} \right.\right\}
$$
and the {\it spectrum} of $A$ by $\sigma(A)={\Bbb C} \backslash \rho(A)$. 
$\lambda$ is an {\it isolated} eigenvalue of  (\ref{var_eig}) if there 
is a neighborhood of $\lambda$, call it $\Omega(\lambda)$, so that
$\Omega(\lambda)\cap \sigma(A)=\{\lambda\}$ (i.e., 
$\lambda$ is an isolated eigenvalue of $A$).  If $\lambda$ is an 
 isolated nonzero eigenvalue of  (\ref{var_eig}) then $Ker[ A - \lambda I]$ 
 is the associated eigenspace.  
${\cal U}=\bigcup_{k=1}^{\infty} Ker[( A - \lambda I)^k]$ similarly 
will be the invariant subspace 
for (\ref{var_eig}) associated with  $\lambda$. No compactness assumptions 
have been made for either $A$ or $T$, so {\it a priori} it may happen 
that (\ref{var_eig}) has {\it no} eigenvalues at all or those that it has may be 
embedded in essential spectrum (defined with respect to $A$) and not isolated. 
 $\lambda$ has finite multiplicity $m$ if $\dim\, {\cal U} =m < \infty$.   
 If $\lambda$ has 
finite multiplicity then there is a finite  integer, $r \leq m $ for which
$Ker[( A - \lambda I)^r]=Ker[( A - \lambda I)^{r+1}]$.  
The smallest such integer is called the {\it ascent } of $\lambda$.

Furthermore, if $\lambda$ is an isolated eigenvalue with finite 
multiplicity, then for each $k=1,\ 2,\ \dots$, $(A-\lambda)^{k}$ is a 
Fredholm operator with zero index implying $nullity(A-\lambda)^{k}=
nullity(A^{*}-\bar{\lambda})^{k}$ for each $k$ (cf. \cite{Kato}, p. 239).  
In particular, $\bar{\lambda}$ is 
an eigenvalue of $A^{*}$ with the same multiplicity and ascent as 
$\lambda$.  A ``left eigenvector'' associated with such a 
$\lambda\in \sigma(A)$ may be characterized variationally as
 \begin{equation}
 \begin{array}{c}
 \mbox{Find  $0\neq \hat{u}^{*}\in Dom(a)$ so that}\\
 	a(\hat{u}^{*}, v)=\lambda \langle \hat{u}^{*}, v \rangle_\scalH 	
 	\mbox{  for all $v\in {Dom(a)}$.}
 	\label{left_eig}
 \end{array}
 \end{equation}
 Note that $ \hat{u}^{*}$  
is a solution for (\ref{left_eig}) if and 
only if $ \hat{u}^{*} \in Dom(A^{*})$ and $\{\bar{\lambda},\ \hat{u}^{*}\}$ is an eigenpair 
for the operator $A^{*}$.  

Henceforth we will assume that there is an isolated 
eigenvalue $\lambda\ne 0$ for (\ref{var_eig}) having finite multiplicity $m$ 
with an associated maximal invariant subspace ${\cal U}$ for which  
we seek approximations. We denote with ${\cal U}^{*}$ the maximal 
invariant subspace for $A^{*}$ associated with $\bar{\lambda}$.

The spectral projection for $A$ onto ${\cal U}$
is defined by the Dunford integral
$$
E=\frac{1}{2\pi\imath}\int_{\Gamma}(z-A)^{-1} dz.
$$ 
where $\Gamma$
is a circle in $\Bbb C$ centered at $\lambda$  
leaving the origin and all  points of $\sigma(A)$ other than $\lambda$ in its exterior. 
The complementary $A$-invariant subspace is denoted ${\cal U}^c=Ran(I-E)$.

Notice that $\mu=1/\lambda$ will be an isolated eigenvalue of $T$ also with multiplicity $m$ 
and the same $m$-dimensional invariant subspace ${\cal U}$ as for $\lambda$.  
The spectral projection may be defined with respect to $T$ as
$$
E=\frac{1}{2\pi\imath}\int_{\Sigma}(z-T)^{-1} dz.
$$ 
where $\Sigma$
is a circle in $\Bbb C$ centered at $\mu$  
leaving the origin and all  points of $\sigma(T)$ other than $\mu$ in its exterior.  

\section{The Galerkin Method}  \label{galerk}
\subsection{Discretization}
In order to approximate  the eigenvalue $\lambda$ and 
its associated invariant subspace ${\cal U}$,
 we introduce two finite dimensional subspaces ${\cal S}_{1,h} \subset {\cal V}$ and  
 ${\cal S}_{2,h}\subset {\cal V}$  -- the trial and test subspaces, respectively.  
 Assume that   $\dim {\cal S}_{1,h}=\dim {\cal S}_{2,h}\stackrel{\rm def}{=}N(h)$.
Typically, the dimension $N(h)$ is monotone increasing as the ``mesh size'' 
parameter $h$ decreases.

The Galerkin method proceeds by solving an eigenvalue problem as in 
(\ref{var_eig}) for the form $a$  restricted to the finite dimensional space ${\cal S}_{2,h}\times {\cal S}_{1,h}$:
 \begin{equation}
 \begin{array}{c}
 \mbox{Find $\lambda_{h}$ and $0\neq v_{h}\in {\cal S}_{1,h}$ so that}\\
 	a(u,v_{h})=\lambda_{h} \langle u, v_{h} \rangle_\scalH \\
 	\mbox{for all $u\in {\cal S}_{2,h}$.}
 \end{array}
 	\label{discvar_eig}
 \end{equation} 
 
 The name is sometimes further qualified as either an 
 orthogonal-Galerkin method or a Petrov-Galerkin method depending on 
 whether ${\cal S}_{1,h}={\cal S}_{2,h}$ or not.  When $A$ is 
 self-adjoint and ${\cal S}_{1,h}={\cal S}_{2,h}$, this approach is 
 called the Rayleigh-Ritz method.
 
  For any given $h$, the computational realization proceeds by fixing
   bases for ${\cal S}_{1,h}$ as 
 $\phi_{1},\ \phi_{2},\ \dots,\ \phi_{N(h)}$, and for ${\cal S}_{2,h}$ as
 $\psi_{1},\ \psi_{2},\ \dots,\ \psi_{N(h)}$.
 The problem (\ref{discvar_eig}) is then reduced to resolving the 
 generalized matrix  eigenvalue problem
\begin{align}
	\bA_{h}\by =&  \lambda_{h} \bB_{h}\by \label{matrixeig}\\
	 \mbox{where  } \bA_{h}&=[a(\psi_{i},\ \phi_{j})]\in {\Bbb C}^{N(h)\times 
	 N(h)} \nonumber \\  
	 \mbox{  and  } \bB_{h}&=
	 [\langle \psi_{i},\ \phi_{j} \rangle_\scalH]\in {\Bbb C}^{N(h)\times N(h)}.
  \nonumber
  \end{align}
 If an eigenvector $\by$ of (\ref{matrixeig}) has components $\by^{t}=\{y_{1},\ 
 y_{2},\ \dots,\ y_{N(h)}\}$, then the corresponding $v_{h}$ that solves 
 (\ref{discvar_eig}) is represented as $v_{h}=\sum_{j=1}^{N(h)}y_{j}\phi_{j}$.
 
 For any $\tau\in \Bbb C$, define
 $ T_{\tau}=\,T+\tau$, which may be defined variationally by analogy to 
 (\ref{BTdef}) as that operator that satisfies
  \begin{equation}
	 \langle u,v\rangle_\scalH+\tau a(u,v)  =a(u,T_{\tau}v)
  	\label{WLOG}
  \end{equation}
 for all $u,\, v\in {\cal  V}$. 
 Notice that $\{\lambda,\ u\}$ is an eigenpair for (\ref{var_eig}) if and 
only if $\{(\lambda^{-1}+\tau),\ u\}$ is an eigenpair for the operator 
$T_{\tau}$ and more generally,  $\sigma(T_{\tau})=  \sigma(T)+\tau$.
 $T$ and $T_{\tau}$ have the same invariant 
subspaces ${\cal U}$ associated with each of the eigenvalues $\lambda^{-1}$ 
and $\lambda^{-1}+\tau$, respectively.  
  The effect of a translation of $\sigma(T)$ by $\tau$
   produces (from (\ref{WLOG}))
 a  discrete problem with translated  spectrum. 
 Instead of (\ref{matrixeig}), we have
 \begin{equation}
	 \bA_{h}\by = {\hat \lambda_{h}}(\bB_{h}+\tau \bA_{h})\by.
 	\label{TtauDisc}
 \end{equation}
 The approximate spectra produced by (\ref{matrixeig}) and (\ref{TtauDisc}) 
 are related as $ {\hat \lambda_{h}}^{-1}=\lambda_{h}^{-1}+\tau$ but 
  eigenvectors and invariant subspaces are identical.  Since our 
 principal interest is in eigenvector approximations, choices for 
 $\tau$ are immaterial, and particular choices will entail no loss of generality.  
 
  Assume  that the following 
 ``discrete inf-sup'' conditions are satisfied:
 \begin{align}
	\inf \begin{Sb}
	       u\in {\cal S}_{2,h}\\
	       \|u\|_\scalV = 1
	       \end{Sb}
	   \sup \begin{Sb}
	       v\in {\cal S}_{1,h} \\
	       \|v\|_\scalV = 1
	       \end{Sb}
    |a(u,v)| & \stackrel{\rm def}{=} \beta(h) >0. 	\label{discinf_supa}\\
	\intertext{and}
	\inf \begin{Sb}
	       u\in {\cal S}_{2,h}\\
	       \|u\|_\scalH = 1
	       \end{Sb}
	   \sup \begin{Sb}
	       v\in {\cal S}_{1,h} \\
	       \|v\|_\scalH = 1
	       \end{Sb}
	    |\langle u,v\rangle_{\cal H} | & \stackrel{\rm def}{=} 
	    \overset{\circ}{\beta}(h) >0. 	\label{discinf_supb}
\end{align}
Since $\dim {\cal S}_{1,h}=\dim {\cal S}_{2,h}=N(h)$, 
	these are equivalent to the complementary conditions,	
	 \begin{align}
	 \inf \begin{Sb}
	       v\in {\cal S}_{1,h}\\
	       \|v\|_\scalV = 1
	       \end{Sb}
	       \sup \begin{Sb}
	       u\in {\cal S}_{2,h}\\
	       \|u\|_\scalV = 1
	       \end{Sb}
	        |a(u,v)| & = \beta(h) >0  \label{discinf_supa2}\\
	\intertext{and}
	\inf \begin{Sb}
	      v\in {\cal S}_{1,h} \\
	       \|v\|_\scalH = 1
	       \end{Sb}
	   \sup \begin{Sb}
	        u\in {\cal S}_{2,h}\\
	       \|u\|_\scalH = 1
	       \end{Sb}
	    |\langle u,v\rangle | & =
	    \overset{\circ}{\beta}(h) >0, \label{discinf_supb2}	
 \end{align}
respectively. Condition (\ref{discinf_supa}) is the usual discrete inf-sup condition
(cf. \cite{BabOsb}) and guarantees that $\bA_{h}$ 
 is invertible for each  $h$.  Analogously condition (\ref{discinf_supb}) 
  guarantees that $\bB_{h}$ 
 is invertible for each  $h$.   Either (\ref{discinf_supa}) or (\ref{discinf_supb})
 will guarantee that the discrete eigenvalue problem 
 (\ref{matrixeig}) is well-posed and associated with a regular matrix 
 pencil for each $h>0$.
 
\subsection{Projections}
 Define  $P_{h}: {\cal  V} \rightarrow {\cal S}_{1,h}$  
 as $P_{h}v=\sum_{i,j=1}^{N(h)}a(\psi_{j},v)\gamma_{ij}\phi_{i}$ 
 where $[\gamma_{ij}]=\bA_{h}^{-1}$. 
  Direct calculation verifies  $P_{h}^{2}=P_{h}$, hence
  $P_{h}$ is  a projection, albeit 
 nonorthogonal typically.  
 $P_{h}$ maps each $v\in {\cal V}$ 
 to a unique vector,  $w^{\sharp}=P_{h}v$ in  ${\cal S}_{1,h}$ that solves
 \begin{equation}
 \begin{array}{c}
 \mbox{Find $w^{\sharp}\in {\cal S}_{1,h}$ so that}\\ 
 a(u,v-w^{\sharp})=0 \quad \mbox{for all $u\in {\cal S}_{2,h}$.}
 	\end{array}
 	\label{Pdef}
 \end{equation}
 $P_{h}$ arises spontaneously in discussing solutions to boundary value 
 problems associated with $a$.
 For any given $f\in {\cal H}$, the weakly posed boundary value problem
 \begin{equation*}
 \begin{array}{c}
 \mbox{Find $\hat{v}\in {\cal V}$ so that}\\
 	a(w,\hat{v})= \langle w, f \rangle_\scalH 	
 	\mbox{  for all $w\in {\cal V}$.}
 \end{array}
 \end{equation*}
 admits a solution $\hat{v}$ which may be approximated with a Galerkin 
 method 
 \begin{equation*}
 \begin{array}{c}
 \mbox{Find $\hat{v}_{h}\in {\cal S}_{1,h}$ so that}\\
 	a(w,\hat{v}_{h})= \langle w, f \rangle_\scalH 	
 	\mbox{  for all $w\in {\cal S}_{2,h}$.}
 \end{array}
 \end{equation*} 
 Exact and approximate solutions are related via $P_{h}$ as $\hat{v}_{h}=P_{h}\hat{v}$.
 
 Along the same lines as above, define $P_{h}^a$ as 
$P_{h}^{a}u=\sum_{i,j=1}^{N(h)}\overline{a(u,\phi_{i})}\overline{\gamma}_{ij}
\psi_{j}$.  $P_{h}^{a}$ is a projection onto  ${\cal S}_{2,h}$ defined on
${\cal V}$ and $w^{\sharp}=P_{h}^{a}u$ solves, for any $u\in {\cal  V}$,
 \begin{equation}
 \begin{array}{c}
 \mbox{Find $w^{\sharp}\in {\cal S}_{2,h}$ so that}\\ 
	 a(u-w^{\sharp},v)=0 \quad \mbox{for all $v\in {\cal S}_{1,h}$.}
 	\end{array}
 	\label{Padef}
 \end{equation}

  Notice that (\ref{Pdef}) and (\ref{Padef}) together imply 
  for all $u,\ v\ \in {\cal  V}$, 
 \begin{equation*}
	 a(u,P_{h}v)=a(P_{h}^{a}u,P_{h}v)=a(P_{h}^{a}u,v)
 \end{equation*}
 That is, $P_{h}^{a}$ is the ``$a$-adjoint''  of $P_{h}$.
 
Now, for all $u\in {\cal S}_{2,h}$ and all $v \in {\cal S}_{1,h}$, 
 \begin{align*}
	 \langle u,\ v \rangle_\scalH=& \langle P_{h}^{a}u,P_{h}v \rangle_\scalH  \\
	           &=a(P_{h}^{a}u,TP_{h}v) \\
	           &=a(u,P_{h}TP_{h}v), 
 \end{align*}
 so we have that  
 $\lambda_{h}\neq 0$ and $v_{h}$ together solve (\ref{discvar_eig})
 if and only if $\lambda_{h}^{-1}$ and $v_{h}$ constitute an eigenpair for
 $T_{h}\stackrel{\rm def}{=}P_{h}TP_{h}$.
 
 From (\ref{discinf_supa2}), 
we find for any $v\in {\cal V}$ with $\|P_{h}v\|_\scalV\neq 
0 $,
\begin{align*}
	  0< \beta(h) \leq & \sup \begin{Sb}
	          u\in {\cal S}_{2,h} \\
	              \|u\|_\scalV = 1
	              \end{Sb}
	    \frac{|a(u,P_{h}v)|}{\|P_{h}v\|_\scalV} = \sup \begin{Sb}
	          u\in {\cal S}_{2,h} \\
	              \|u\|_\scalV = 1
	              \end{Sb}
	    \frac{|a(u,v)|}{\|P_{h}v\|_\scalV} \\
	    \leq & \sup \begin{Sb}
	          u\in {\cal S}_{2,h} \\
	              \|u\|_\scalV = 1
	              \end{Sb} \frac{c_{1}\|v\|_\scalV\|u\|_\scalV}
	              {\|P_{h}v\|_\scalV}    \leq    c_{1}	\frac{\|v\|_\scalV}{\|P_{h}v\|_\scalV}
\end{align*}
Thus, 
\begin{equation}
	\|P_{h}\|_\scalV \leq c_{1}/\beta(h) .
	\label{projbnd}
\end{equation}
Similarly from (\ref{discinf_supb}), 
 $\|P_{h}^{a}\|_\scalV \leq c_{1}/\beta(h) $.
 The following result leads us to conclude that both
 $\|I-P_{h}\|_\scalV \leq c_{1}/\beta(h) $ and  
 $\|I- P_{h}^{a}\|_\scalV \leq c_{1}/\beta(h) $ as well.
\begin{lemma} \label{idemp}
If $Z$ is a bounded (nonorthogonal) projection on a Hilbert space 
${\cal W}$, then $\|I-Z\|_{\cal W}=\|Z\|_{\cal W}$.  Furthermore,
if $\Pi$ denotes the ${\cal W}$-orthogonal projection onto $Ran(Z)$,
\begin{equation}
	\frac{1}{\|Z\|_{\cal W}}\|(I-Z)u\|_{\cal W}\leq \|(I-\Pi)u\|_{\cal W} 
	\leq \|(I-Z)u\|_{\cal W}
	\label{Zbound}
\end{equation}
\end{lemma}
{\sc Proof:}
The first assertion was proved by Kato (\cite{Kato0}, p. 28).  Since
$(I-Z)=(I-Z)(I-\Pi)$, $\|(I-Z)u\|_{\cal W}\leq \|(I-Z)\|_{\cal W}\|(I-\Pi)u\|_{\cal W}$
which then gives the first inequality of (\ref{Zbound}).  The second 
inequality of (\ref{Zbound}) is the best approximation property of 
orthogonal projections. $\quad 	  \blacksquare $
 
$\Pi_{1,h}$ and $\Pi_{2,h}$ will always denote orthogonal projections onto
${\cal S}_{1,h}$ and ${\cal S}_{2,h}$ respectively. However,  depending on
the context, they will be considered either orthogonal in ${\cal H}$ or
orthogonal in ${\cal V}$ with no distinction in notation.

 Define  $Q_{h}: {\cal  H} \rightarrow {\cal S}_{1,h}$  
 as $Q_{h}v=\sum_{i,j=1}^{N}\langle \psi_{j},v\rangle_{\cal 
 H}\, \overset{\circ}{\gamma}_{ij}\phi_{i}$ 
 where $[\, \overset{\circ}{\gamma}_{ij}]=\bB_{h}^{-1}$. 
 $Q_{h}$ has a natural extension to $v\in {\cal  V}'$ so the composition of 
 operators $AQ_{h}: {\cal  V}' \rightarrow {\cal V}'$ and 
 $Q_{h}AQ_{h}: {\cal  V}' \rightarrow {\cal S}_{1,h}$ are each well defined.
 Since $Q_{h}^{2}=Q_{h}$, $Q_{h}$ is  also a projection, but is ${\cal
H}$-orthogonal if and only if ${\cal S}_{1,h}={\cal S}_{2,h}$.  
 $Q_{h}$ maps each $v\in {\cal H}$ 
 to a unique vector,  $w^{\sharp}=Q_{h}v$ in  ${\cal S}_{1,h}$ that solves
 \begin{equation}
 \begin{array}{c}
 \mbox{Find $w^{\sharp}\in {\cal S}_{1,h}$ so that}\\ 
	 \langle u,v-w^{\sharp}\rangle_\scalH=0 \quad \mbox{for all $u\in {\cal S}_{2,h}$.}
 	\end{array}
 	\label{Qdef}
 \end{equation}
 Evidently, the ${\cal H}$-adjoint $Q_{h}^*: {\cal  H}\rightarrow {\cal 
 S}_{2,h}$ has the form 
$Q_{h}^{*}u=\sum_{i,j=1}^{N} \langle \phi_{i},u\rangle_{\cal 
 H}\, \overline{\overset{\circ}{\gamma}_{ij}}\psi_{j}$.  $Q_{h}^{*}$ is a projection onto  
 ${\cal S}_{2,h}$ 
and solves, for any $u\in {\cal H}$,
 \begin{equation}
 \begin{array}{c}
 \mbox{Find $w^{\sharp}\in {\cal S}_{2,h}$ so that}\\ 
	 \langle u-w^{\sharp},v\rangle_\scalH=0 \quad \mbox{for all $v\in {\cal S}_{1,h}$.}
 	\end{array}
 	\label{Qstardef}
 \end{equation}

Now, for all $u\in {\cal S}_{2,h}$ and all $v \in {\cal S}_{1,h}$, 
 \begin{align*}
	 a(u,v)=& a( Q_{h}^{*}u,Q_{h}v)  \\
	           &=\langle Q_{h}^{*}u,AQ_{h}v\rangle_\scalH \\
	           &=\langle u,Q_{h}AQ_{h}v\rangle_\scalH, 
 \end{align*}
 so we have that
 $\lambda_{h}\neq 0$ and $v_{h}$ together solve (\ref{discvar_eig})
 if and only if $\lambda_{h}$ and $v_{h}$ constitute an eigenpair for
 $A_{h}\stackrel{\rm def}{=}Q_{h}AQ_{h}$.

\subsection{Convergence}
Convergence criteria may be framed either in ${\cal V}$ or in  ${\cal H}$.
Convergence criteria in  ${\cal V}$ appear as 
\begin{align}
	\lim_{h\rightarrow 0} \beta(h)^{-1} 
	\inf_{w\in {\cal S}_{1,h}} & \| v-w\|_\scalV=0 \quad 
	\mbox{for each $v \in {\cal V}$,}
	\label{V_S1approx}\\
\intertext{and}
		\lim_{h\rightarrow 0} 
      \delta_\scalV(T_{h}, T) =0,
	       \label{T_S1approx}
\end{align}

\begin{theorem}[ Descloux, {\it et al.} \cite{Des1978,Des1981}] 
\label{DNR}
The hypotheses (\ref{V_S1approx})  and (\ref{T_S1approx}) imply that
\begin{enumerate}
	\item Both $P_{h}\rightarrow I$ and $P_{h}^{a}\rightarrow I$ strongly 
in ${\cal V}$; $P_{h}$ is uniformly ${\cal V}$-bounded  with respect to $h$; and
there is a constant $c>0$ so that
$$
\delta_\scalV(T_{h}, T) \leq \|(I-P_{h})TP_{h}\|_\scalV\leq 
c\, \delta_\scalV(T_{h}, T).
$$  
   \item  For each compact subset, ${\cal R}$, of $\rho(T)$ there exists
	$h_{0}>0$ and $K>0$ so that ${\cal R}\subset \rho(T_{h})$ 
    and $\|(z-T_{h})^{-1}\|_\scalV<K$ uniformly for 
	$z \in {\cal R}$ for all $h<h_{0}$.
	\item If $\mu$ is an eigenvalue of $T$ with algebraic multiplicity $m$ 
	and with an associated $m$-dimensional invariant 
	subspace ${\cal U}$, there will be $m$ eigenvalues (counting 
	multiplicity) of $T_{h}$, $\{ \mu_h^{1},\ 
	\mu_h^{2},\ \dots,\ \mu_h^{m}\} $ that are convergent to $\mu$ as $h\rightarrow 0$ 
	and the associated $m$-dimensional $T_{h}$-invariant subspace ${\cal U}_{h}$ 
	satisfies $\delta({\cal U},{\cal U}_{h})\rightarrow 0$ as $h\rightarrow 0$.
\end{enumerate}
\end{theorem}

If $\Sigma$
is a circle in $\Bbb C$ centered at $\mu$  
leaving the origin and all  points of $\sigma(T)$ other than $\mu$ in its 
exterior, then under the convergence
assumptions (\ref{V_S1approx}) and (\ref{T_S1approx}),  there will be 
$m$ eigenvalues of $T_{h}$, labeled as $\mu_h^{1},\ \mu_h^{2},\ \dots,\ 
\mu_h^{m}, $ that will all be  
contained in the interior of $\Sigma$ for sufficiently small $h$.   
  Thus for sufficiently small $h$, the Dunford integral
$$
E_{h}=\frac{1}{2\pi\imath}\int_{\Sigma}({\hat z}-T_h)^{-1} d{\hat z}
$$
defines a spectral projection onto the $T_h$-invariant subspace 
${\cal U}_{h}$ associated with $\mu_h^{1},\ \mu_h^{2},\ \dots\ \mu_h^{m}$.

 Analogous convergence criteria in ${\cal H}$ appear as
 \begin{align}
	\lim_{h\rightarrow 0} \overset{\circ}{\beta}(h)^{-1} 
	\inf_{w\in {\cal S}_{1,h}} & \| v-w\|_\scalH=0 \quad 
	\mbox{for each $v \in {\cal H}$,}
	\label{H_S1approx}\\
\intertext{and}
		\lim_{h\rightarrow 0} 
      \delta_\scalH(A_{h}, A) =0,
	       \label{A_S1approx}
\end{align}
with similar consequences:

\begin{theorem}[ Descloux, {\it et al.} \cite{Des1978,Des1981}] 
\label{DLR}
The hypotheses (\ref{H_S1approx})  and (\ref{A_S1approx}) imply that
\begin{enumerate}
	\item Both $Q_{h}\rightarrow I$ and $Q_{h}^{*}\rightarrow I$ strongly 
in ${\cal H}$; $Q_{h}$ is uniformly ${\cal H}$-bounded  with respect to 
$h$; and there is a constant $c>0$ so that
$$
\delta_\scalH(A_{h}, A) \leq \|(I-Q_{h})AQ_{h}\|_\scalH\leq 
c\, \delta_\scalH(A_{h}, A), 
$$ (if $A$ is unbounded the second inequality holds trivially).

\item  For each compact subset ${\cal R}$ of $\rho(A)$, there exists
	$h_{0}>0$ and $K>0$ so that ${\cal R}\subset \rho(A_{h})$ 
    and $\|(z-A_{h})^{-1}\|_\scalH<K$ uniformly for 
	$z \in {\cal R}$ for all $h<h_{0}$.

	\item If $\lambda$ is an eigenvalue of $A$ with algebraic multiplicity $m$ 
	and with an associated $m$-dimensional invariant 
	subspace ${\cal U}$, there will be $m$ eigenvalues (counting 
	multiplicity) of $A_{h}$, $\{ \lambda_h^{1},\ 
	\lambda_h^{2},\ \dots,\ \lambda_h^{m}\} $ that are convergent to $\lambda$ as $h\rightarrow 0$ 
	and the associated $m$-dimensional $A_{h}$-invariant subspace ${\cal U}_{h}$ 
	satisfies $\delta_{\cal H}({\cal U},{\cal U}_{h})\rightarrow 0$ as $h\rightarrow 0$.

\end{enumerate}
\end{theorem}
 
If $\Gamma$ is a circle in $\Bbb C$ centered at $\lambda$  
leaving the origin and all  points of $\sigma(A)$ other than $\lambda$ 
in its exterior, then under the convergence
assumptions (\ref{H_S1approx}) and (\ref{A_S1approx}),  there will be 
$m$ eigenvalues of $A_{h}$, labeled as $\lambda_h^{1},\ \lambda_h^{2},\ \dots,\ 
\lambda_h^{m}, $ that will all be  
contained in the interior of $\Gamma$ for sufficiently small $h$.   
Thus for sufficiently small $h$, the Dunford integral
$$
E_{h}=\frac{1}{2\pi\imath}\int_{\Gamma}({\hat z}-A_h)^{-1} d{\hat z}
$$
defines a spectral projection for $A_h$ onto the invariant subspace 
${\cal U}_{h}$ associated with $\lambda_h^{1},\ \lambda_h^{2},\ \dots\ \lambda_h^{m}$. 

It will be convenient to label the complementary nonzero part of the spectrum of 
$T_{h}$ and $A_{h}$ as $\sigma^c(T_{h})=
\sigma(T_h)\backslash \{0,\  \mu_h^{1},\ \mu_h^{2},\ \dots,\ \mu_h^{m}\}$ 
and $\sigma^c(A_{h})=
\sigma(A_h)\backslash \{0,\  \lambda_h^{1},\ \lambda_h^{2},\ \dots,\ 
\lambda_h^{m}\}$, respectively.
   
\begin{theorem}[Babu\v{s}ka and Osborn \cite{BabOsb}]  
	Suppose the convergence hypotheses (\ref{H_S1approx}) and 
	(\ref{A_S1approx}) hold.  Then 
	\begin{equation} \label{eivecrate}
		\delta_\scalH({\cal U},{\cal U}_{h}) \leq c 
		\left[\overset{\circ}{\beta}(h)^{-1}\delta_\scalH({\cal U},{\cal 
		S}_{1,h})\right]^{1/\alpha}
	\end{equation}
	\begin{equation}
		|\lambda - \lambda_{h}|\leq  c \left[\overset{\circ}{\beta}(h)^{-1}
		\delta_\scalH({\cal U^{*}},{\cal S}_{2,h})\,
		\delta_\scalH({\cal U},{\cal S}_{1,h})\right]^{1/\alpha}
		\label{eigrate}
	\end{equation}
\end{theorem}

These appear as Theorems 8.3 and 8.4 in \cite{BabOsb}. 
As stated there, the proofs given in \cite{BabOsb} presume that $A$ is compact 
 - however, the arguments extend without change to 
unbounded $A$ having  nontrivial essential spectrum once the convergence 
results of \cite{Des1981} come into play.

\section{Bounded $A$} \label{A_bnded}
\subsection{Basic Estimates}
Define 
$$
\varepsilon_\scalH(h)=
\sup_{u\in {\cal U}}\frac{\|Q_{h}A(I-Q_{h})u\|_\scalH}{\|(I-\Pi_{1,h})u\|_\scalH}
$$
where $\Pi_{1,h}$ here is the ${\cal H}$-orthogonal projection onto
${\cal S}_{1,h}$.
\begin{theorem} \label{MainH}
	Suppose the convergence hypotheses (\ref{H_S1approx}) and 
	(\ref{A_S1approx}) hold.  There exists an $h_{0}>0$ sufficiently small 
so that for each $h<h_{0}$ and each $u\in {\cal U}$, there is a 
$u_{h} \in {\cal U}_{h}$ so that
\begin{equation}
	\|u_{h}-Q_{h}u\|_\scalH \leq c\, \varepsilon_\scalH(h)\,
	\delta_\scalH({\cal U},{\cal S}_{1,h})
	\label{projeqn2}
\end{equation}
where $c>0$ is a constant independent of $h$ and independent of the 
choice of $u \in {\cal U}$.
\end{theorem}
{\sc Proof:}    $Ran(E_h)\subset {\cal S}_{1,h}$ 
since $E_hA_h=A_hE_h$.
Note also that $Q_{h}$ is a spectral projection for $A_{h}$ 
associated with all nonzero eigenvalues of $A_{h}$.  Thus, 
$Q_{h}-E_{h}$ is a spectral projection for $A_{h}$ onto ${\cal 
U}^{c}_{h}$ associated with all nonzero eigenvalues of $A_{h}$ 
distinct from $\lambda$.
Let ${\hat A}=A|_{{\cal U}}$ denote the restriction of $A$ to ${\cal U}$ and 
Let ${\hat A}^{c}_{h}=A_{h}|_{{\cal U}^{c}_{h}}$ denote the restriction of $A_{h}$ 
to ${\cal U}^{c}_{h}$.  Then, ${\hat A}^{c}_{h}(Q_{h}-E_{h})=
(Q_{h}-E_{h})A_{h}$ and we have 
\begin{align*}
	{\hat A}^{c}_{h}(Q_{h}-E_{h})\left|_{\scalU}\right. - (Q_{h}-E_{h})\left|_{\scalU}\right.{\hat A} = & 
	(Q_{h}-E_{h})(A_{h}-A)\left|_{\scalU}\right. \\
	=- (Q_{h}-E_{h})&((I-Q_{h})A+Q_{h}A(I-Q_{h}))\left|_{\scalU}\right. \\
	=- & (Q_{h}-E_{h})Q_{h}A(I-Q_{h})\left|_{\scalU}\right. \\
	=- & (I-E_{h})Q_{h}A(I-Q_{h})\left|_{\scalU}\right.
\end{align*}
Thus, the mapping $S:{\cal U} \rightarrow {\cal U}^{c}_{h}$ given by
 $S=(Q_{h}-E_{h})\left|_{\scalU}\right.$ is a solution to the 
 Sylvester equation,
	\begin{equation}
		{\hat A}^{c}_{h}S - S{\hat A} = -(I-E_{h})Q_{h}A(I-Q_{h})\left|_{\scalU}\right.
		\label{syleqnbndA}
	\end{equation}
There exists $K_{1}>0$ such that 
$$
\|(z-{\hat A})^{-1}\left|_{{\cal U}}\right.\|_\scalH\leq
\|(z-A)^{-1}\|_\scalH \leq K_{1}
$$
uniformly for all $z\in \Gamma$.
Likewise there exists an $h_{0}>0$ and $K_{2}>0$ such that for $h<h_{0}$,
$$
\|(z-{\hat A}_{h}^{c})^{-1}\left|_{{\cal U}^{c}_{h}}\right.\|_\scalH\leq
\|(z-A_{h})^{-1}\|_\scalH \leq K_{2}
$$
uniformly for $z\in \Gamma$.  Therefore, the pseudospectral sets 
$\Lambda_{\epsilon}({\hat A}_{h})$ are contained in the exterior of 
$\Gamma$ for any $\epsilon<1/K$ and for all $h>h_{0}$.  By Lemma 
\ref{sepbnd}(b), there must then be a $c_{0}>0$ independent of $h$, such 
that 
$$
\|(Q_{h}-E_{h})\left|_{\scalU}\right.\|_{{\cal U} \rightarrow {\cal 
U}^{c}_{h}} \leq c_{0} \, \|(I-E_{h})Q_{h}A(I-Q_{h})\left|_{\cal 
U}\right.\|_{{\cal U} \rightarrow {\cal U}^{c}_{h}}.
$$
Thus, for any $u\in {\cal U}$,
\begin{align*}
	\|(Q_{h}-E_{h})u\|_\scalH \leq & c_{0}\, \|I-E_{h}\|_\scalH \, 
	\sup_{w \in {\cal U}}\frac{\|Q_{h}A(I-Q_{h})w\|_{\cal 
	V}}{\|w\|_\scalH} \\
    \leq &	c_{0}\, \|E_{h}\|_\scalH \, 
	\sup_{w \in {\cal U}}\frac{\|Q_{h}A(I-Q_{h})w\|_{\cal 
	V}}{\|(I-\Pi_{1,h})w\|_\scalH}\, 
	\sup_{w \in {\cal U}} \frac{\|(I-\Pi_{1,h})w\|_{\cal	V}}
	{\|w\|_\scalH}\\
	= & c_{0}\, \|E_{h}\|_\scalH \, \varepsilon_\scalH(h)\,
	\delta_\scalH({\cal U},{\cal S}_{1,h}).
\end{align*}
Notice that since $E_{h}$ converges uniformly to $E$, $\|E_{h}\|_\scalH $ is 
uniformly bounded.  The conclusion follows upon assigning 
$u_{h}=E_{h}u$.  $\quad \blacksquare$

 \begin{corollary} \label{sacase_Abnded}
  Suppose the convergence hypotheses (\ref{H_S1approx}) and 
	(\ref{A_S1approx}) hold and that ${\cal S}_{1,h}={\cal S}_{2,h}$.
 \begin{equation}
 	\delta_\scalH({\cal U},{\cal S}_{1,h}) \leq 
 	\delta_\scalH({\cal U},{\cal U}_{h})
	\leq (1+c\,  \varepsilon_\scalH(h) )\delta_\scalH({\cal U},{\cal S}_{1,h})
		\label{close_Abnded}
 \end{equation}
\end{corollary}

{\sc Proof:}
Note that under the hypotheses given, $Q_{h}=\Pi_{1,h}=\Pi_{2,h} 
=Q_{h}^{*}$.
The first inequality of (\ref{close_Abnded}) follows trivially 
from observing that ${\cal U}_{h}\subset{\cal S}_{1,h}$.
For the second, by Theorem \ref{MainH} there exists an $h_{0}>0$ such that for each 
$h<h_{0}$ and $u\in {\cal U}$ with $\|u\|_\scalH=1$, there is a ${\hat u}_{h}\in {\cal U}_{h}$
such that $\|Q_{h}u-{\hat u}_{h}\|_\scalH\leq c\, \varepsilon_\scalH(h)\,
	\delta_\scalH({\cal U},{\cal S}_{1,h})$.  Then,
	\begin{align*}
		\min_{u_{h} \in {\cal U}_{h}} \|u-u_{h}\|_\scalH 
		\leq & \|u-{\hat u}_{h}\|_\scalH \\
		\leq & \|u-Q_{h}u\|_\scalH+\|Q_{h}u-{\hat u}_{h}\|_\scalH \\
		\leq & \delta_\scalH({\cal U},{\cal S}_{1,h})+c\, \varepsilon_\scalH(h)\,
		\delta_\scalH({\cal U},{\cal S}_{1,h})
			\end{align*}
Maximizing over $u$ yields the conclusion. $ \quad \blacksquare$

Corollary \ref{sacase_Abnded} can be interpretted as saying that, provided 
$\varepsilon_\scalH(h) \rightarrow 0$, the orthogonal Galerkin method 
provides optimal zero order approximations to the eigenvectors of $A$:
 ${\cal U}_{h}$ will approach 
asymptotically the {\it closest} $m$-dimensional subspace in ${\cal S}_{1,h}$ to the exact 
eigenspace ${\cal U}$ -- and this is true even if $A$ is nonselfadjoint. 
See the comments following Theorem \ref{sacase}.

\subsection{Related Estimates and Interpretation}
	Since each $Q_{h}$ is a bounded projection, ${\cal H}$ may be 
decomposed into a direct sum of complementary subspaces as ${\cal H}=
Ran(Q_{h})\oplus Ker(Q_{h})={\cal S}_{1,h}\oplus{\cal S}_{2,h}^{\perp}$.
The operator $A$ can then be partitioned in a way that reflects this 
decomposition -- see Figure 2.

For each generalized eigenvector, $u \in {\cal U}$, define $u^{(1,h)}=Q_hu$
and $u^{(2,h)}=(I-Q_h)u$.
The quantity 
$\varepsilon_\scalH(h)$ then is a measure of the relative size of $Q_{h}A(I-Q_{h})$ 
on $u^{(2,h)}$.  This follows from 
observing that from Lemma \ref{idemp} 
$$
\sup_{u\in{\cal U}}\frac{\|Q_{h}A(I-Q_{h})u\|)}{\|(I-Q_{h})u\|}\leq 
\varepsilon_\scalH(h)\leq \|Q_{h}\|\sup_{u\in{\cal U}}\frac{\|Q_{h}A(I-Q_{h})u\|)}{\|(I-Q_{h})u\|},
$$
which (provided $Q_{h}\rightarrow I$ strongly) implies immediately
$$
\sup_{u^{(2,h)}\in{\cal N}_{h}({\cal 
U})}\frac{\|Q_{h}A(I-Q_{h})u^{(2,h)}\|_\scalH}{\|u^{(2,h)}\|_\scalH}\leq 
\varepsilon_\scalH(h)\leq c\, \sup_{u^{(2,h)}\in{\cal N}_{h}({\cal 
U})}\frac{\|Q_{h}A(I-Q_{h})u^{(2,h)}\|_\scalH}{\|u^{(2,h)}\|_\scalH}$$
for some constant $c>0$, where 
${\cal N}_{h}({\cal U})=
(I-Q_{h}){\cal U}\subset {\cal S}_{2,h}^{\perp}$.  ${\cal N}_{h}({\cal
U})$ is the span of all components $u^{(2,h)}$ of (generalized)
eigenvectors $u$ in ${\cal U}$ that lie the direction of $Ker(Q_{h})$.

\begin{center}
\setlength{\unitlength}{2.7mm}
\begin{picture}(40,30)(-5,-5)
\thicklines
\put(-1,-1){\line(0,1){22}}
\put(21,-1){\line(0,1){22}}
\put(-1,-1){\line(1,0){1.5}}
\put(21,-1){\line(-1,0){1.5}}
\put(-1,21){\line(1,0){1.5}}
\put(21,21){\line(-1,0){1.5}}
\put(28,-1){\line(0,1){22}}
\put(32,-1){\line(0,1){22}}
\put(28,-1){\line(1,0){.75}}
\put(32,-1){\line(-1,0){.75}}
\put(28,21){\line(1,0){.75}}
\put(32,21){\line(-1,0){.75}}
\thinlines
\put(19.9,15){\makebox(0,0){\psboxit{box .9 setgray fill}
    {\hspace*{5.3cm}\rule[-2.65cm]{0cm}{5.3cm}}}}
\put(0,0){\dashbox(10,20){}}
\put(0,10){\dashbox(20,10){}}
\put(0,0){\dashbox(20,10){}}
\put(-3,10){\makebox(0,0){$A\ =$}}
\put(5,15){\makebox(0,0){$A_{h}$}}
\put(15,15){\makebox(0,0){\scriptsize{$Q_{h}A(I-Q_{h})$}}}
\put(5,5){\makebox(0,0){\scriptsize{$(I-Q_{h})AQ_{h}$}}}
\put(15,5){\makebox(0,0){\scriptsize{$(I-Q_{h})A(I-Q_{h})$}}}
\put(31.35,5){\makebox(0,0){\psboxit{box .9 setgray fill}
    {\hspace*{1.49cm}\rule[-2.65cm]{0cm}{5.3cm}}}}
\put(28.5,0){\dashbox(3,10){}}
\put(28.5,10){\dashbox(3,10){}}
\put(25.5,10){\makebox(0,0){$u\ =$}}
\put(30.2,15){\makebox(0,0){\scriptsize{$u^{(1,h)}$}}}
\put(30.2,5){\makebox(0,0){\scriptsize{$u^{(2,h)}$}}}

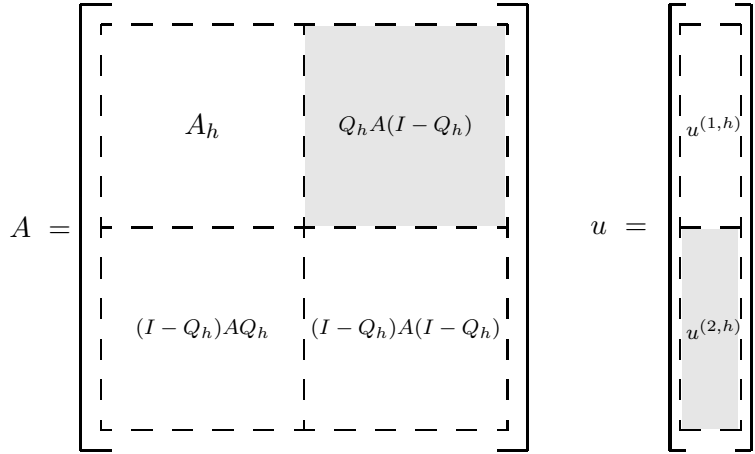
\put(-4,-6){\shortstack{Figure 2.  Partitioning $A$ and $u\in{\cal U}$
on ${\cal H}={\cal S}_{1,h}\oplus {\cal S}_{2,h}^{\perp}$  \\
Shading indicates influence on $\varepsilon_\scalH(h)$}}
\end{picture}  
\end{center}
\vspace{2mm}

Certainly it may happen that $\varepsilon_\scalH(h)\not\rightarrow 0$, so a
variety of additional conditions will be examined in the next few sections that
suffice to guarantee $\varepsilon_\scalH(h)\rightarrow 0$. Perhaps the simplest of
these is to require that $A_h^*$ eventually converge to $A^*$ in gap:
 \begin{theorem} \label{epsbnd_Abnded}
 There is a $c>0$ such that
$
	 \varepsilon_\scalH(h) \leq  c \delta_\scalH(A_{h}^{*},A^{*})
 $
 \end{theorem}
 {\sc Proof:}
Note that  $I-Q_{h}=(I-Q_{h})(I-\Pi_{1,h})$.
Thus, 
\begin{align*}
\varepsilon_\scalH(h) & = \sup_{u\in {\cal U}}\frac{\|Q_{h}A(I-Q_{h})u\|_\scalH}
{\|(I-\Pi_{1,h})u\|_\scalH}    \\
&  \leq \|Q_{h}A(I-Q_{h})\|_\scalH  
 = \|(I-Q_{h}^{*})A^{*}Q_{h}^{*}\|_\scalH \\
	           & \leq  c \delta_\scalH(A_{h}^{*},A^{*}), 
\end{align*}
using Part 1 of Theorem \ref{DLR}. $\quad \blacksquare$

We shouldn't expect to do much better than the bound provided by (\ref{projeqn2}).
 The bound has the ``right'' asymptotic behaviour in many cases and so in 
that sense will be best possible. 

\begin{theorem} \label{Necessary}
Suppose that $\lambda$ is a simple eigenvalue of $A$ with an associated
eigenvector $u$. There exist constants
$c_0,\, c_{1}>0$ (independent of $h$) so that for each $h$ one may find
a Galerkin eigenvector $u_h$ with
\begin{equation}
c_0 \|u_{h}-Q_{h}u\|_\scalH \leq \|Q_{h}A(I-Q_{h})u\|_\scalH \leq c_{1}\|u_{h}-Q_{h}u\|_\scalH 
+|\lambda - \lambda_{h}|\, \|u_{h}\|_\scalH
	\label{nec_eqn1}
\end{equation}
Furthermore there exists a
$c_{2}>0$ independent of $h$ so that
\begin{equation}
	c_0 \frac{\|u_{h}-Q_{h}u\|_\scalH}
	{\|(I-\Pi_{1,h})u\|_\scalH}\leq \varepsilon_\scalH(h) \leq c_{1} \frac{\|u_{h}-Q_{h}u\|_\scalH}
	{\|(I-\Pi_{1,h})u\|_\scalH}+(c_{2}/_{\overset{\circ}{\beta}(h)})
	\|(I-\Pi_{2,h})u^{*}\|_\scalH
	\label{nec_eqn2}
\end{equation}
where $u^{*}$ is a ``left eigenvector'' satisfying (\ref{left_eig}).
\end{theorem}
 {\sc Proof:}
The first inequalities in each of (\ref{nec_eqn1}) and (\ref{nec_eqn2}) are a
consequence of Theorem \ref{MainH}.
Since $\lambda$ is simple, $Au=\lambda u$ and for $h>0$ sufficiently 
small,  $rank(E_{h})=1$ and $\lambda_{h}$ will be simple.
 \begin{align*}
	 \|Q_{h}A(I-Q_{h})u\|_\scalH \leq & \|(I-E_{h})Q_{h}A(I-Q_{h})u\|_\scalH + 
	 \|E_{h}Q_{h}A(I-Q_{h})u\|_\scalH \\
	  \leq & \|({\hat A}^{c}_{h}S - S{\hat A})u\|_\scalH + 
	 \|E_{h}Au-E_{h}A_{h}u\|_\scalH \\
	 \leq & \left(\|{\hat A}^{c}_{h}\|_\scalH + |\lambda | \right) \|Su\|_\scalH + 
	 \|\lambda E_{h}u-A_{h}E_{h}u\|_\scalH \\
	 \leq & \left(\| A_{h}\|_\scalH + |\lambda | \right) \|Su\|_\scalH + 
	 \|(\lambda -\lambda_{h})E_{h}u\|_\scalH \\
	 \leq & \| A\|_\scalH\left( 1 + \| Q_{h}\|_\scalH^{2} \right) \|E_{h}u-Q_{h}u\|_\scalH + 
	 |\lambda -\lambda_{h}| \|E_{h}u\|_\scalH,
 \end{align*}
 where $S=(Q_{h}-E_{h})\left|_{\scalU}\right.$ satisfies the Sylvester equation (\ref{syleqnbndA}).
(\ref{nec_eqn1}) then follows upon assigning $u_{h}=E_{h}u$ and 
observing that $Q_{h}$ is uniformly bounded in ${\cal H}$.

To show (\ref{nec_eqn2}), note that the ascent of $ \lambda $ is 1, so 
from (\ref{eigrate}) there is a $c_{2}$ with
$$
		|\lambda - \lambda_{h}|\leq  c_{2} \overset{\circ}{\beta}(h)^{-1}
		\|(I-\Pi_{2,h})u^{*}\|_\scalH
		\|(I-\Pi_{h})u\|_\scalH
$$
so
\begin{align*}
	\varepsilon_\scalH(h) = & \frac{\|Q_{h}A(I-Q_{h})u\|_\scalH}
		{\|(I-\Pi_{h})u\|_\scalH} \\
	\leq & c_{1} \frac{\|u_{h}-Q_{h}u\|_\scalH}
		{\|(I-\Pi_{h})u\|_\scalH} + (c_{2}/_{\overset{\circ}{\beta}(h)})
			\|(I-\Pi_{2,h})u^{*}\|_\scalH \quad \blacksquare
\end{align*}

\section{Unbounded $A$ -- Estimates in ${\cal V}$} \label{Vest}
\subsection{Basic results}
The setting considered in this section is the traditional one 
encountered in the analysis of finite element methods.  With  few exceptions,
much of the structure of arguments of Section 4 carry over into this setting.

 Define 
$$
\varepsilon_\scalV(h)=
\sup_{u\in {\cal U}}\frac{\|P_{h}T(I-P_{h})u\|_{\cal 
V}}{\|(I-\Pi_{1,h})u\|_\scalV}
$$
where $\Pi_{1,h}$ now is the ${\cal V}$-orthogonal projection onto
${\cal S}_{1,h}$.

\begin{theorem} \label{MainV}
Suppose the convergence hypotheses 
(\ref{V_S1approx})-(\ref{T_S1approx}) hold.  There exists an $h_{0}>0$ sufficiently small 
so that for each $h<h_{0}$ and all $u\in {\cal U}$, there is a 
$u_{h} \in {\cal U}_{h}$ so that
\begin{equation}
	\|u_{h}-P_{h}u\|_\scalV \leq c\, \varepsilon_\scalV(h)\,
	\delta_\scalV({\cal U},{\cal S}_{1,h})
	\label{projeqn3}
\end{equation}
where $c>0$ is a constant independent of $h$ and independent of the 
choice of $u \in {\cal U}$.
\end{theorem}
The proof is the same as for Theorem \ref{MainH} formulated in the 
Hilbert space ${\cal V}$ instead of ${\cal H}$ and with $T$, $T_{h}$, 
and $P_{h}$ playing the roles of $A$, $A_{h}$, 
and $Q_{h}$, respectively.

\begin{lemma}   \label{Tcompact} Suppose the convergence hypotheses 
(\ref{V_S1approx})-(\ref{T_S1approx}) hold.  If ${\cal   V}$ is 
compactly imbedded in ${\cal   H}$ (so that $T$ is compact as a 
mapping from ${\cal V}$ to itself) and if ${\cal S}_{2,h}$ 
satisfies the approximation property
	\begin{equation}
		\lim_{h\rightarrow 0} \beta(h)^{-1} 
		\inf_{w\in {\cal S}_{2,h}} \| v-w\|_\scalV=0 \quad 
		\mbox{for each } v \in {\cal V},
		\label{S2approx}
	\end{equation}
then $\varepsilon_\scalV(h) \rightarrow 0$ as $h\rightarrow 0$.
\end{lemma}
{\sc Proof:} Since $a$ is bounded and coercive on ${\cal V}$, there is a 
bounded and invertible
linear operator on ${\cal V}$, ${\cal A}$, such that $a(u,v)=
\langle u, {\cal A}v \rangle_\scalV$. Let ``${\star}$'' denote the ${\cal V}$-adjoint
 and observe that 
$P_{h}^{\star}u=\sum_{i,j=1}^{N}\langle \phi_{i},u \rangle_\scalV 
\overline{\gamma_{ij}} {\cal A}^{\star}\psi_{j}$, so that, in 
particular, $Ran(P_{h}^{\star})={\cal A}^{\star}{\cal S}_{2,h}$. 
Let $\Pi_{2,h}$ denote the ${\cal V}$-orthogonal 
projection onto ${\cal S}_{2,h}$ and notice that ${\cal A}^{\star}\Pi_{2,h}{\cal A}^{-\star}$ 
is a projection (no longer orthogonal, in general) 
onto $Ran(P_{h}^{\star})$.  Then for any $u\in {\cal V}$, (\ref{S2approx}) implies 
\begin{align*}
	\|(I-P_{h}^{\star})u\|_\scalV & =\|(I-P_{h}^{\star})(I-{\cal A}^{\star}\Pi_{2,h}{\cal A}^{-\star})u\|_\scalV \\
	       & =\|(I-P_{h}^{\star}){\cal A}^{\star}(I-\Pi_{2,h}){\cal A}^{-\star}u\|_\scalV \\
	        & \leq \|I-P_{h}^{\star}\|_\scalV\|{\cal A}^{\star}\|_\scalV\|(I-\Pi_{2,h}){\cal A}^{-\star}u\|_\scalV \\
           & \leq \frac{c}{\beta(h)}\|(I-\Pi_{2,h}){\cal A}^{-\star}u\|_\scalV \\
         & \quad =	\frac{c}{\beta(h)}
\inf_{w\in {\cal S}_{2,h}}  \| {\cal A}^{-\star} u-w\|_\scalV\rightarrow 0, 
\end{align*}
for some constant $c$.
Thus, $P_{h}^{\star}$ converges strongly to $I$ in ${\cal V}$.  
Since $T^{\star}$ is compact, $\|(I-P_{h}^{\star})T^{\star}\|_\scalV\rightarrow 0$ and
\begin{align*}
	\varepsilon_\scalV(h) & = \sup_{u\in {\cal U}}\frac{\|P_{h}T(I-P_{h})(I-\Pi_{1,h})u\|_{\cal 
V}}{\|(I-\Pi_{1,h})u\|_\scalV}    \\
&  \leq \|P_{h}T(I-P_h)\|_\scalV  \\
& = \|(I-P_{h}^{\star})T^{\star}P_{h}^{\star}\|_\scalV \\
	           & \leq \|(I-P_{h}^{\star})T^{\star}\|_\scalV 
	           \|P_{h}^{\star}\|_\scalV\rightarrow 0 
	           \quad 	  \blacksquare 
\end{align*}

Even when $T$ is not compact, additional conditions on ${\cal S}_{2,h}$ 
can yield the same result:
\begin{lemma}  \label{Tgen}
 Let $T^{*}$ denote the ${\cal H}$-adjoint of $T$.  Suppose the convergence hypotheses 
(\ref{V_S1approx})-(\ref{T_S1approx}) hold.  If  ${\cal S}_{2,h}$ 
satisfies the approximation properties
\begin{align}
	\lim_{h\rightarrow 0} \beta(h)^{-1} 
		\inf_{w\in {\cal S}_{2,h}} \| v-w\|_\scalV & =0 \quad 
		\mbox{for each } v \in {\cal V},
		\label{S2approxb} \\
\intertext{and}
			\sup \begin{Sb}
		       v\in {\cal S}_{2,h} \\
		       \|v\|_\scalV = 1
		       \end{Sb}  \inf_{w\in {\cal S}_{2,h}}  \| T^{*}v-w\|_\scalV 
		       \stackrel{\rm def}{=} \gamma(h) & \rightarrow 0 \quad \mbox{ as $h\rightarrow 0$}.
		       \label{T_S2approxb}
\end{align}
Then   $\varepsilon_\scalV(h) \rightarrow 0$ as $h\rightarrow 0$ and 
$\varepsilon_\scalV(h) = {\cal O}(\gamma(h))$.
\end{lemma}
{\sc Proof:}   We first verify that $T^{*}$ maps ${\cal V}$ to ${\cal V}$.
For any $u\in {\cal H},\ v\in Dom(a)$,  
$|\langle u,\, v\rangle_\scalH|\leq \frac{1}{\alpha}\|u\|_\scalH\ 
\|C_{a}v\|_\scalH$.
Thus,  for $u\in Dom(a),\ w\in Ran(C_{a}|_{Dom(A)})$, 
(so that $w=C_{a}v$ for some $v\in 
Dom(A)$), one may observe
\begin{align*}
	|\langle T^{*}u,\, B_{a}^{*}w\rangle_\scalH|& = |\langle T^{*}u,\, B_{a}^{*}C_{a}v
	\rangle_\scalH| \\
	&= |\langle T^{*}u,\, Av\rangle_\scalH| \\
	&=|\langle u,\, v\rangle_\scalH|  \\
	 & \leq \left(\frac{\|u\|_\scalH}{\alpha}\right) 
	\|w\|_\scalH
\end{align*}
Since $Dom(A)$ is a core for $C_{a}$, $Ran(C_{a}|_{Dom(A)})$ is dense 
in ${\cal H}$ and as a consequence $T^{*}u\in 
Dom(B^{**}_{a})=Dom(B_{a})=Dom(a).$
For $u\in Dom(a),\ v\in Dom(A)$, then:
\begin{align}
	a(T^{a}u,v)=a(u,Tv)=&\langle u,\, v\rangle_\scalH\\ \nonumber
	=&\langle u,\, TAv\rangle_\scalH\\  \nonumber
	=&\langle T^{*}u,\, Av\rangle_\scalH=a(T^{*}u,v), \label{a_adj}
\end{align}
(the last equality being a consequence of $T^{*}u\in Dom(a)$) and so, $T^{a}=T^{*}$.

Since $\|\cdot\|_\scalV $ and $[\Re e\,a]$ are equivalent 
norms on ${\cal V}$, there is an $m>0$ so that $m\|u\|_\scalV^{2}
\leq |a(u,u)|$ and
\begin{align*}
m \|P_{h}T(I-P_{h})v\|_\scalV & \leq \frac{|a(P_{h}T(I-P_{h})v,\ 
            P_{h}T(I-P_{h})v)|}{\|P_{h}T(I-P_{h})v\|_\scalV}\\
   & \leq \sup\begin{Sb}
		       u\in {\cal V} \\
		       \|u\|_\scalV = 1
		       \end{Sb} |a(u,\ P_{h}T(I-P_{h})v)| \\
		       &  = \sup\begin{Sb}
		       u\in {\cal V} \\
		       \|u\|_\scalV = 1
		       \end{Sb} |a((I-P_{h}^{a})T^{*}P_{h}^{a}u,\ v)| \\
		       &\leq \sup\begin{Sb}
		       u\in {\cal V} \\
		       \|u\|_\scalV = 1
		       \end{Sb} c_{1}\|v\|_\scalV\|(I-P_{h}^{a})T^{*}P_{h}^{a}u\|_\scalV \\
		       &\quad = c_{1}\|v\|_\scalV\|(I-P_{h}^{a})T^{*}P_{h}^{a}\|_\scalV
\end{align*}
Thus, $\varepsilon_\scalV(h)\leq \|P_{h}T(I-P_{h})\|_\scalV\leq 
(c_{1}/m)\|(I-P_{h}^{a})T^{*}P_{h}^{a}\|_\scalV$.  Now, notice that 
$Ran(P_{h}^{a})={\cal S}_{2,h}$, so $I-P_{h}^{a}=(I-P_{h}^{a})(I-\Pi_{2,h})$
-- where $\Pi_{2,h}$ is the ${\cal V}$-orthogonal projection onto 
${\cal S}_{2,h}$.
Now,
\begin{align*}
\|(I-P_{h}^{a})T^{*}P_{h}^{a}\|_\scalV & \leq 
\|(I-P_{h}^{a})(I-\Pi_{2,h})T^{*}\Pi_{2,h}P_{h}^{a}\|_\scalV \\
 & \leq \|(I-P_{h}^{a})\|_\scalV\|P_{h}^{a}\|_\scalV\|(I-\Pi_{2,h})T^{*}
 \Pi_{2,h}\|_\scalV\\
 &\quad = \|(I-P_{h}^{a})\|_\scalV\|P_{h}^{a}\|_\scalV \sup \begin{Sb}
		       v\in {\cal S}_{2,h} \\
		       \|v\|_\scalV = 1
		       \end{Sb}  \inf_{w\in {\cal S}_{2,h}}  \| T^{*}v-w\|_\scalV \\
		       &\quad \leq c \gamma(h) \rightarrow 0  \quad \blacksquare 
\end{align*}

If $a$ is symmetric (so that $[\Re e\ a]=a$ and $[\Im m\ a]=0$)
and if $a(\cdot,\ \cdot)$ itself is used as the inner 
product for ${\cal V}$, then  $T$ is a self-adjoint operator in ${\cal V}$.
If additionally ${\cal S}_{1,h}={\cal S}_{2,h}$, then $P_{h}$ is a
${\cal V}$-orthogonal projection and $T_{h}$ is then also 
 self-adjoint.  In this circumstance, $u_{h}$ is asymptotically 
 the {\em closest} vector out of ${\cal S}_{1,h}$ to $u$ (with respect 
 to the $a$-norm on ${\cal V}$):
 \begin{theorem} \label{sacase}
 Suppose $a$ is symmetric,  $a(u,v)=\langle u,\, v\rangle_\scalV$, and
  that ${\cal S}_{1,h}={\cal S}_{2,h}$.
Then
 \begin{equation}
	 \varepsilon_\scalV(h) \leq  c_{1} \delta_\scalV(T_{h},T)
 	\label{epsbnd}
 \end{equation}
and
 \begin{equation}
 	1 \leq 
 	\frac{\delta_\scalV({\cal U},{\cal U}_{h})}{\delta_\scalV({\cal U},{\cal S}_{1,h})}
	\leq 1+c_{2}\, \delta_\scalV(T_{h},T)
		\label{close}
 \end{equation}
\end{theorem}
{\sc Proof:}
Note that under the hypotheses given, $P_{h}=\Pi_{1,h}=\Pi_{2,h}=P_{h}^{a}$.
Thus, 
\begin{align*}
	\varepsilon_\scalV(h) & = \sup_{u\in {\cal U}}\frac{\|\Pi_{1,h}T(I-\Pi_{1,h})u\|_{\cal 
V}}{\|(I-\Pi_{1,h})u\|_\scalV}    \\
&  \leq \|\Pi_{1,h}T(I-\Pi_{1,h})\|_\scalV  \\
& = \|(I-\Pi_{1,h})T\Pi_{1,h}\|_\scalV \\
	           & \leq  c_{1} \delta_\scalV(T_{h},T)
\end{align*}
The first inequality of (\ref{close}) follows trivially 
from observing that ${\cal U}_{h}\subset{\cal S}_{1,h}$.
For the second, by Theorem \ref{MainV} there exists an $h_{0}>0$ such that for each 
$h<h_{0}$ and $u\in {\cal U}$ with $\|u\|_\scalV=1$, there is a ${\hat u}_{h}\in {\cal U}_{h}$
such that $\|P_{h}u-{\hat u}_{h}\|_\scalV\leq c\, \varepsilon_\scalV(h)\,
	\delta_\scalV({\cal U},{\cal S}_{1,h})$.  Then,
	\begin{align*}
		\min_{u_{h} \in {\cal U}_{h}} \|u-u_{h}\|_\scalV 
		\leq & \|u-{\hat u}_{h}\|_\scalV \\
		\leq & \|u-P_{h}u\|_\scalV+\|P_{h}u-{\hat u}_{h}\|_\scalV \\
		\leq & \delta_\scalV({\cal U},{\cal S}_{1,h})+c\, \varepsilon_\scalV(h)\,
		\delta_\scalV({\cal U},{\cal S}_{1,h})
			\end{align*}
Maximizing over $u$ yields the conclusion. $ \quad \blacksquare$

Theorem \ref{sacase} was essentially given by  Chatelin \cite{Chat},  refined 
by Babu\v{s}ka and Osborn \cite{BabOsb} -- each for compact self-adjoint $T$. 
Recently, a more general result of this sort allowing for noncompact self-adjoint $T$ was 
given by Knyazev \cite{Kny}.

\subsection{Elliptic boundary value problems: Finite Elements}
 Let $\Omega$ be a bounded open subset of 
 ${\Bbb R}^{n}$ with  a
 boundary $\partial \Omega$ that is at least $C^{r+1,1}$ for some integer $r>0$.  
 Given real coefficient functions 
$a_{ij},\  b_{i},\ c \in C^{r}(\bar{\Omega})$, consider the second order 
elliptic differential operator $A$ defined  by 

\begin{align*}
	A({\bf x},D)v= & -\sum_{i,\ j=1}^{n}D_{j}a_{ij}({\bf x})D_{i}v 
	        + \sum_{i=1}^{n}b_{i}({\bf x})D_{i}v + c({\bf x})v \quad 
	        \mbox{in  $\Omega$} 
\end{align*}
with $v\,=\,0$  on $\partial \Omega$,
 and the related adjoint operator
given by
\begin{align*}
A^{*}({\bf x},D)u= & -\sum_{i,\ j=1}^{n}D_{i}a_{ij}({\bf x})D_{j}u 
        - \sum_{i=1}^{n}D_{i}(b_{i}({\bf x})u) + c({\bf x})u  \quad 
	        \mbox{in  $\Omega$} 
\end{align*}
with $u\,=\,0$  on $\partial \Omega$,

Suppose that $A({\bf x},D)$ is uniformly strongly elliptic.
 The associated bilinear form:
\begin{align*}
a(w,\ v)=\sum_{i,\ j=1}^{n}\int_{\Omega}a_{ij}({\bf x})
[D_{j}w] & [D_{i}v] d{\bf x}\ + \sum_{j=1}^{n}\int_{\Omega} w\ b_{j}({\bf x})
[D_{j}v] d{\bf x} \\
& +  \int_{\Omega} c({\bf x})\ w\ v\, d{\bf x},
\end{align*}
defined on $Dom(a)=H_{0}^{1}(\Omega)$
is a closed sectorial bilinear form densely defined in ${\cal H}=L^{2}(\Omega)$ and 
$A$ is manifested as a densely defined m-sectorial operator on ${\cal H}$
which can be extended to
a continuous bijection from ${\cal V}= H_{0}^{1}(\Omega)$ 
onto the dual space ${\cal V}'=H^{-1}(\Omega)$.  Here and elsewhere, 
$H^{p}(\Omega)$ denotes the completion of the vector space 
$C^{\infty}$ with respect to the norm 
$$
\|u\|_{H^{p}(\Omega)}=\sum_{|\alpha|\leq 
p}\int_{\Omega}|D^{\alpha}v|^{2}d{\bf x}.
$$
  The associated 
seminorm is defined as 
$$
\nrm\, u \,\nrm_{\mbox{\raisebox{-.4ex}{\scriptsize 
$H^{p}(\Omega)$}}}=\sum_{|\alpha|= 
p}\int_{\Omega}|D^{\alpha}v|^{2}d{\bf x}.
$$

Results governing regularity of solutions to elliptic problems
(e.g., \cite{wloka} p. 328) guarantee for any $f\in H^{r-1}(\Omega)$, 
the weakly posed problem
\begin{align}
	a(w,v)&=\langle w,f\rangle_{L^{2}(\Omega)}\quad 
	\mbox{for all $w \in H_{0}^{1}(\Omega)$} 
\end{align}
has a solution $v\in H^{r+1}(\Omega)\cap H^{1}_{0}(\Omega)$.
Since $T$ is a closed mapping from $H^{r-1}(\Omega)$ to $H^{r+1}(\Omega)$,
there exists a constant $c>0$
 such that $\|Tf\|_{H^{r+1}(\Omega)}\leq c \|f\|_{H^{r-1}(\Omega)}$
 for all $f\in H^{r-1}(\Omega)$.
Furthermore, if ${\cal U}$ denotes an invariant subspace of $A$ associated with 
an isolated eigenvalue of $A$ with finite multiplicity then 
${\cal U}\subset H^{r+1}(\Omega)$.  

Likewise the adjoint problem,
	\begin{equation}
		a(u,w)=\langle w,g\rangle_{L^{2}(\Omega)}\quad 
		\mbox{for all $w \in H^{1}(\Omega)$} 
		\label{adjoint}
	\end{equation}
has a solution $u\in H^{r+1}(\Omega)$ for any $g\in H^{r-1}(\Omega)$.
So, in particular, if $g\in {\cal V}\subset H^{1}(\Omega)$, then 
$T^{*}g=(A^{*})^{-1}g\in H^{3}(\Omega)$.
 
Apply the Galerkin method with ${\cal S}_{h}={\cal S}_{1,h}=
{\cal S}_{2,h}$ chosen to be a family of finite dimensional subspaces of ${\cal V}$, 
 so that for all  integers $0<k \leq r$ and some fixed $c>0$, $u\in 
 H^{r+1}(\Omega)$ implies
	\begin{equation}
		\inf_{v\in {\cal S}_{h}}\|u-v\|_{H^{1}(\Omega)}\leq c\, h^{k}\nrm\, u 
		\,\nrm_{\mbox{\raisebox{-.4ex}{\scriptsize $H^{k+1}(\Omega)$}}}
		\label{interp}
	\end{equation}
For example, $C^{0}$-finite element spaces constructed from piecewise 
polynomials of degree at least $r$ would satisfy this condition.

The discrete inf-sup condition (\ref{discinf_supa}) is satisfied 
with $\beta(h)=\alpha>0$. Thus, the convergence condition \ref{V_S1approx}
is immediately satisfied. It remains to verify that (\ref{T_S1approx}) 
holds.  Note that for every $x\in {\cal S}_{h}=Dom(T_{h})$,
$Tx\in H^{3}(\Omega)$ and
     \begin{align*}
      	 \inf_{\,y\in H^{1}(\Omega)} 
		   \frac{\| x -  y \|_{H^{1}(\Omega)}+\| T_{h}x -  Ty \|_{H^{1}(\Omega)} }
		   { \|x\|_{H^{1}(\Omega)}+\|T_{h}x\|_{H^{1}(\Omega)}} &\leq  
		   \frac{\| T_{h}x -  Tx \|_{H^{1}(\Omega)} }
		   { \|x\|_{H^{1}(\Omega)}+\|T_{h}x\|_{H^{1}(\Omega)}} \\
		  & = \frac{\| (I-P_{h})Tx \|_{H^{1}(\Omega)} }
		   { \|x\|_{H^{1}(\Omega)}+\|T_{h}x\|_{H^{1}(\Omega)}} \\
		   & \leq \frac{c\, h^{2} \nrm\, Tx \,\nrm{\mbox{\raisebox{-.6ex}{\scriptsize $H^{3}(\Omega)$}}} }
		   { \|x\|_{H^{1}(\Omega)}+\|T_{h}x\|_{H^{1}(\Omega)}} \\
		   & \leq \frac{c\, h^{2} \| x \|_{H^{1}(\Omega)} }
		   { \|x\|_{H^{1}(\Omega)}+\|T_{h}x\|_{H^{1}(\Omega)}} \\
		   & \leq c\ h^{2}
      \end{align*}
Thus, $\lim_{h\rightarrow 0} 
      \delta_\scalV(T_{h}, T) =0$.
 
Although Theorem \ref{Tcompact} is applicable, Theorem \ref{Tgen} will yield a concrete rate 
once we evaluate
$$
\gamma(h)=\sup \begin{Sb}
		       v\in {\cal S}_{h} \\
		       \|v\|_\scalV = 1
		       \end{Sb}  \inf_{w\in {\cal S}_{h}}  \| T^{*}v-w\|_\scalV 		      
$$
For any $v\in {\cal S}_{h} $, $T^{*}v \in H^{3}(\Omega)$ and so
\begin{align*}
	\inf_{w\in {\cal S}_{h}}  \| T^{*}v-w\|_\scalV & =\|(I-\Pi_{h}) T^{*}v\|_\scalV \\
	      & \leq c h^{2} \nrm\,T^{*}v \,\nrm{\mbox{\raisebox{-.6ex}{\scriptsize $H^{3}(\Omega)$}}} \\
	      & \leq c h^{2} \|v\|_{H^{1}(\Omega)}
\end{align*}
Thus, $\gamma(h)\leq c h^{2}$.

Since ${\cal U}\subset H^{r+1}(\Omega)$, 
      \begin{align*}
	      \delta_{\cal  V}({\cal  U},{\cal S}_{h}) &  =  
		\sup_{ x\in{\cal  U}}\inf_{\,y\in{\cal S}_{h}} 
		   \frac{\| y -  x \|_\scalV}{\| x\|_\scalV} 
		   =\frac{\|(I-\Pi_{h}) x\|_\scalV}{\| x\|_\scalV} \\
	      & \leq c h^{r} \frac{\nrm\, x \,\nrm{\mbox{\raisebox{-.6ex}
	      {\scriptsize $H^{r+1}(\Omega)$}}} }{\| x\|_\scalV} \leq c h^{r}.
      \end{align*}

Theorem \ref{MainV} asserts that there exists an $h_{0}>0$ sufficiently small 
so that for each $h<h_{0}$ and all $u\in {\cal U}$, there is a 
$u_{h} \in {\cal U}_{h}$ so that
$$
	\|u_{h}-P_{h}u\|_\scalV \leq c\, h^{r+2},
	$$
whereas $\|u_{h}-u\|_\scalV$ and $\|u-P_{h}u\|_\scalV$ will each 
be only of order $h^{r}$ typically.

\section{Unbounded $A$ -- Estimates in ${\cal H}$} \label{Hest}
\subsection{Basic results}
In the ${\cal V}$-setting explored in Section 
5, orthogonality of $P_{h}$ and the related best approximation 
property in ${\cal V}$  could be developed only for self-adjoint $A$.
In contrast, estimates in ${\cal H}$ such as were found in Section 4 
have particular appeal since whenever ${\cal S}_{1,h}={\cal S}_{2,h}$, $Q_{h}$ will  be an
{\em orthogonal} projection in  ${\cal H}$, notwithstanding asymmetry in $a$
and nonselfadjointness of $A$. 
Unfortunately, those estimates obtained  in Section 4 depend 
fundamentally on the boundedness of $A$. In particular, if $A$ is unbounded then
$\varepsilon_\scalH(h)$ might not be uniformly bounded in $h$, much less go
to zero. 

We define an expression that plays a role analogous to that of 
$\varepsilon_\scalH(h)$ in Section 4: 
\begin{equation}
	\overset{\circ}{\varepsilon}_\scalH(h)=
	\sup_{u\in {\cal U}}\frac{\|(P_{h}-Q_{h})u\|_\scalH}
	{\|(I-\Pi_{1,h})u\|_\scalH}
	\label{epsH}
\end{equation}
where $\Pi_{1,h}$ here is once again the ${\cal H}$-orthogonal projection onto
${\cal S}_{1,h}$.

\begin{theorem} \label{MainH_unb}
	Suppose the convergence hypotheses (\ref{H_S1approx}) and 
	(\ref{A_S1approx}) hold. There exists an $h_{0}>0$ sufficiently small 
so that for each $h<h_{0}$ and all $u\in {\cal U}$, there is a 
$u_{h} \in {\cal U}_{h}$ so that
\begin{equation}
	\|u_{h}-Q_{h}u\|_\scalH \leq c\, \overset{\circ}{\varepsilon}_\scalH(h)\,
	\delta_\scalH({\cal U},{\cal S}_{1,h})
	\label{projeqnH}
\end{equation}
where $c>0$ is a constant independent of $h$ and independent of the 
choice of $u \in {\cal U}$.
\end{theorem}
{\sc Proof:}    
Since $E_hA_h=A_hE_h$, $Ran(E_h)\subset {\cal S}_{1,h}$
Note also that $Q_{h}$ is a spectral projection for $A_{h}$ 
associated with all nonzero eigenvalues of $A_{h}$.  Thus, 
$Q_{h}-E_{h}$ is a spectral projection for $A_{h}$ onto ${\cal U}^{c}_{h}$.
Let ${\hat A}=A|_{{\cal U}}$ denote the restriction of $A$ to ${\cal U}$ and 
Let ${\hat A}^{c}_{h}=A_{h}|_{{\cal U}^{c}_{h}}$ denote the restriction of $A_{h}$ 
to ${\cal U}^{c}_{h}$.  Then, ${\hat A}^{c}_{h}(Q_{h}-E_{h})=
(Q_{h}-E_{h})A_{h}$ and we have 
\begin{align*}
	{\hat A}^{c}_{h}(Q_{h}-E_{h})\left|_{\scalU}\right. - 
	(Q_{h}-E_{h})\left|_{\scalU}\right.{\hat A} = & 
	(Q_{h}-E_{h})(A_{h}-A)\left|_{\scalU}\right.  \\
	=- & (Q_{h}-E_{h})((I-Q_{h})A+Q_{h}A(I-Q_{h}))\left|_{\scalU}\right.  \\
	=- & (Q_{h}-E_{h})Q_{h}A(I-Q_{h})\left|_{\scalU}\right.  \\
	=- & (I-E_{h})Q_{h}A(I-Q_{h})\left|_{\scalU}\right.
\end{align*}
Now premultiply by $T_{h}$ and postmultiply by $T$ to find
\begin{align*}
	(Q_{h}-E_{h})\left|_{\scalU}\right.{\hat T} - 
	{\hat T}^{c}_{h}(Q_{h}-E_{h})\left|_{\scalU}\right. = 
	& - T_{h}(I-E_{h})Q_{h}A(I-Q_{h})T\left|_{\scalU}\right. \\
	= & - (I-E_{h})T_{h}Q_{h}A(I-Q_{h})T\left|_{\scalU}\right. \\
	= & - (I-E_{h})(T_{h}Q_{h}AT-T_{h}Q_{h}AQ_{h}T)\left|_{\scalU}\right. \\
	= & - (I-E_{h})(T_{h}Q_{h}-Q_{h}T)\left|_{\scalU}\right. \\
	= & - (I-E_{h})(P_{h}T-Q_{h}T)\left|_{\scalU}\right.
\end{align*}
Thus, the mapping $S:{\cal U} \rightarrow {\cal U}^{c}_{h}$ given by
 $S=(Q_{h}-E_{h})\left|_{\scalU}\right.$ is a solution to the 
 Sylvester equation,
$$
	{\hat T}^{c}_{h}S - S{\hat T} = (I-E_{h})(P_{h}-Q_{h})T\left|_{\scalU}\right.
$$

The goal now is to show that the bounds developed in Appendix A are 
applicable.  There exists a $K_{1}>0$ such that 
$$
\|(z-{\hat A})^{-1}\left|_{{\cal U}}\right.\|_\scalH\leq
\|(z-A)^{-1}\|_\scalH \leq K_{1}
$$
uniformly for all $z\in \Gamma$.
Likewise there exists an $h_{0}>0$ and $K_{2}>0$ such that for $h<h_{0}$,
$$
\|(z-{\hat A}_{h}^{c})^{-1}\left|_{{\cal U}^{c}_{h}}\right.\|_\scalH\leq
\|(z-A_{h})^{-1}\|_\scalH \leq K_{2}
$$
uniformly for $z\in \Gamma$.  Therefore, the pseudospectral sets 
$\Lambda_{\epsilon}({\hat A}_{h})$ are contained in the exterior of 
$\Gamma$ for any $\epsilon<1/K_{2}$ and for all $h>h_{0}$.  By Lemma 
\ref{sepbnd}(b), there must then be a $c>0$ independent of $h$, such 
that 
$$
\|(Q_{h}-E_{h})\left|_{\scalU}\right.\|_{{\cal U} \rightarrow {\cal 
U}^{c}_{h}} \leq c \, \|(I-E_{h})Q_{h}A(I-Q_{h})\left|_{\cal 
U}\right.\|_{{\cal U} \rightarrow {\cal U}^{c}_{h}}.
$$
Thus, for any $u\in {\cal U}$,
\begin{align*}
	\|(Q_{h}-E_{h})u\|_\scalH \leq & c\, \|I-E_{h}\|_\scalH \, 
	\sup_{w \in {\cal U}}\frac{\|Q_{h}A(I-Q_{h})w\|_\scalH}{\|w\|_\scalH} \\
    \leq &	c\, \|E_{h}\|_\scalH \, 
	\sup_{w \in {\cal U}}\frac{\|Q_{h}A(I-Q_{h})w\|_\scalH}{\|(I-\Pi_{1,h})w\|_\scalH}\, 
	\sup_{w \in {\cal U}} \frac{\|(I-\Pi_{1,h})w\|_{\cal	V}}
	{\|w\|_\scalH}\\
	= & c\, \|E_{h}\|_\scalH \, \overset{\circ}{\varepsilon}_\scalH(h)\,
	\delta_\scalH({\cal U},{\cal S}_{1,h}).
\end{align*}
Since $E_{h}$ converges uniformly to $E$, $\|E_{h}\|_\scalH $ is 
uniformly bounded.  The conclusion follows upon assigning 
$u_{h}=E_{h}u$.  $\quad \blacksquare$

When $A$ is unbounded on ${\cal H}$ (so that ${\cal V}\neq{\cal H}$),
$\|P_{h}\|_\scalH$ will not typically be uniformly bounded with 
respect to $h$.  Estimating the rate at which 
$\overset{\circ}{\varepsilon}_\scalH(h)\rightarrow 0$ as $h\rightarrow 0$ becomes 
technically more demanding and additional hypotheses are honorably
acquired.  For the remainder of the section we assume 
that  
$$
{\cal S}_{1,h}={\cal S}_{2,h}\stackrel{\rm def} {=} 
{\cal S}_{h} =span\{ \phi_{1},\ \phi_{2},\ \dots,\ \phi_{N(h)}\}
$$ and  that 
$$
\mbox{{\em The trial vectors $\phi_{i}$ are eigenvectors of 
an auxiliary operator $A_{0}$}}.
$$
 \begin{theorem} \label{Hcase}
 Suppose that  $A$ is decomposable as $A=A_{0}+B$ so that $A_{0}$ 
 is positive definite and self-adjoint on $Dom(A_{0})=Dom(A)$ and $B^{*}$ 
 is bounded relative to $A_{0}$ with relative $A_{0}$-bound less than 1. 
 Furthermore, assume that $\lambda$ is a nondefective eigenvalue of $A$. 
 \begin{enumerate}
 \item  
 If $B^{*}$ is relatively compact with respect to $A_{0}$ 
 then  $\overset{\circ}{\varepsilon}_\scalH(h)\rightarrow 0$ as $h\rightarrow 0$. 
 \label{RelCmpt}
 \item If 
 $$
 \sup_{v\in {\cal S}_{h}} \inf_{w \in {\cal S}_{h}} 
  \frac{\|B^{*}v-w\|_\scalH}{\|A_{0}v\|_\scalH} \stackrel{\rm def} {=}
  \overset{\circ}{\gamma}(h)\rightarrow 0
  $$
 then $\overset{\circ}{\varepsilon}_\scalH(h)\rightarrow 0$ as $h\rightarrow 0$ and 
 $\overset{\circ}{\varepsilon}_\scalH(h) \leq {\cal O}(\overset{\circ}{\gamma}(h))$ as $h\rightarrow 0$.
 \label{RelBnded}
 \end{enumerate}
\end{theorem}

{\sc Proof:}
Since $B^{*}$ has relative bound with respect to $A_{0}$ smaller than 1, there is a $\tau \in 
{\Bbb C}$ so that $\|B^{*}(A_{0}-\tau)^{-1}\|_\scalH=1-\kappa<1$ 
for some $\kappa>0$.
It will be useful to translate the spectrum of $A$ by $\tau$ and write
$A-\tau=(A_{0}-\tau) + B$.  Referring to the discussion around 
(\ref{WLOG}), we absorb this shift in spectrum into both $A$ and 
$A_{0}$, and assume without loss of generality 
that $\|B^{*}A_{0}^{-1}\|_\scalH=1-\kappa<1$.  Before continuing, we 
verify that the assumptions of (\ref{RelCmpt}) and (\ref{RelBnded}) 
are preserved for any appropriate choice of $\tau$.
First considering (\ref{RelCmpt}), if $\tau$ is in the resolvent set for $A_{0}$ then 
$\tau^{-1}$ is in the resolvent set for $A_{0}^{-1}$ and
$$
B^{*}(A_{0}-\tau)^{-1}=\frac{1}{\tau}B^{*}A_{0}^{-1}(\tau^{-1}-A_{0}^{-1})^{-1}.
$$
$B^{*}A_{0}^{-1}$ is compact, $(\tau^{-1}-A_{0}^{-1})^{-1}$ is 
bounded, so $B^{*}(A_{0}-\tau)^{-1}$ is also compact and the 
assumption of (\ref{RelCmpt}) will be independent of any feasible 
shift $\tau$.
Likewise for (\ref{RelBnded}), there is an $m>0$ so that
for any $v\in Dom(A_{0})$, $\|(A_{0}-\tau)v\|_\scalH\geq 
m\|v\|_\scalH$. Thus for any $v\in Dom(A_{0})$,
$$
\|A_{0}v\|_\scalH\leq \|(A_{0}-\tau)v\|_\scalH+|\tau|\cdot\|v\|_\scalH
\leq (1+\frac{|\tau|}{m})\|(A_{0}-\tau)v\|_\scalH.
$$
If the assumption of (\ref{RelBnded}) holds, then
$$
\sup_{v\in {\cal S}_{h}} \inf_{w \in {\cal S}_{h}} 
  \frac{\|B^{*}v-w\|_\scalH}{\|(A_{0}-\tau)v\|_\scalH} \leq (1+\frac{|\tau|}{m})
  \overset{\circ}{\gamma}(h)\rightarrow 0
  $$
so the assumption of (\ref{RelBnded}) is independent of the selected 
shift $\tau$.

Now we prove (\ref{RelBnded}) first.
Write $A=A_{0}(I+A_{0}^{-1}B)$ and observe then that
\begin{align*}
	T_{h}Q_{h}A(I-Q_{h})Tu= & P_{h}TQ_{h}A(I-Q_{h})Tu \\
	                      = & P_{h}TA_{0}Q_{h}A_{0}^{-1}A(I-Q_{h})Tu \\
	                      = & P_{h}TA_{0}Q_{h}A_{0}^{-1}B(I-Q_{h})Tu
\end{align*}
Note that 
\begin{align*}
	\|P_{h}TA_{0}v\|_\scalH^{2}
	      = & \|\sum_{i,j=1}^{N(h)} a(\phi_{j},TA_{0}v)\gamma_{ij}\phi_{i}\|_\scalH^{2} \\
	      = & \sum_{i=1}^{N(h)} |\sum_{j=1}^{N(h)}\langle \phi_{j}, v\rangle_\scalH 
	      \lambda_{j}^{0}\gamma_{ij}|^{2} \\
	            = & \|{\hat \bG}\bx\|_{{\Bbb C}^{N(h)}}^{2}
\end{align*}
where 
$$
\bx=\{\langle \phi_{1}, v\rangle_\scalH,\ \langle \phi_{2}, 
v\rangle_\scalH,\ \dots,\ \langle \phi_{N(h)}, v\rangle_\scalH\}^{t}
$$ 
and ${\hat \bG}$ is the matrix inverse to 
\begin{align*}
	[(\lambda_{\ell}^{0})^{-1}a(\phi_{\ell}, \phi_{k})] = & 
	[\langle A_{0}^{-1}\phi_{\ell},A\phi_{k}\rangle_\scalH] \\
	= & [\langle \phi_{\ell},(I+ A_{0}^{-1}B)\phi_{k}\rangle_\scalH] \\
	= & {\bf I} + [\langle B^{*}A_{0}^{-1}\phi_{\ell}, \phi_{k}\rangle_\scalH]. 
\end{align*}
Observe that for any $\by\in {\Bbb C}^{N(h)}$ with $\|\by\|_{{\Bbb C}^{N(h)}}=1$,
$$
\|[\langle B^{*}A_{0}^{-1}\phi_{\ell}, \phi_{k}\rangle_\scalH]
\by\|_{{\Bbb C}^{N(h)}}\leq \|B^{*}A_{0}^{-1}\|_\scalH\leq 1- \kappa,
$$
hence  we obtain a bound for $P_{h}TA_{0}$ that is uniform in $h$:
$$
\|P_{h}TA_{0}\|_\scalH\leq \|{\hat \bG}\|_{{\Bbb C}^{N(h)}}\leq 
\frac{1}{1-(1-\kappa)}=\frac{1}{\kappa}.
$$
Since $\lambda$ is nondefective, $\lambda Tu=u$ for any $u\in {\cal U}$ 
and we find
\begin{align*}
	\overset{\circ}{\varepsilon}_\scalH(h) = & 
	\sup_{u\in {\cal U}}\frac{\|T_{h}Q_{h}A(I-Q_{h})Tu\|_\scalH}
	{\|(I-Q_{h})u\|_\scalH} \\
	\leq & \sup_{u\in {\cal U}}
	\frac{\|Q_{h}A_{0}^{-1}B(I-Q_{h})Tu\|_\scalH}
	{\kappa\|(I-Q_{h})u\|_\scalH} \\
	\leq & \sup_{u\in {\cal U}}
	\frac{\|Q_{h}A_{0}^{-1}B(I-Q_{h})\cdot (I-Q_{h})Tu\|_\scalH}
	{\kappa|\lambda|\cdot\|(I-Q_{h})Tu\|_\scalH}  \\
		\leq & 	\frac{1}{\kappa|\lambda|} \|Q_{h}A_{0}^{-1}B(I-Q_{h})\|_\scalH.
\end{align*}
Now if the condition of (\ref{RelCmpt}) above holds then 
\begin{align*}
	\|Q_{h}A_{0}^{-1}B(I-Q_{h})\|_\scalH 
	    \leq & \|A_{0}^{-1}B(I-Q_{h})\|_\scalH \\
	     = & \|(I-Q_{h})B^{*}A_{0}^{-1}\|_\scalH.
\end{align*}
Since $B^{*}A_{0}^{-1}$ is compact and $Q_{h}\rightarrow I$ strongly,
$\|(I-Q_{h})B^{*}A_{0}^{-1}\|_\scalH\rightarrow 0$ as $h\rightarrow 0$
from which follows the conclusion of (\ref{RelCmpt}).
     
Now observe that $A_{0}^{-1}{\cal S}_{h}={\cal S}_{h}$ so that, if the 
assumption (\ref{RelBnded}) holds, then
\begin{align*}
	\|Q_{h}A_{0}^{-1}B(I-Q_{h})\|_\scalH
	     = & \|(I-Q_{h})B^{*}A_{0}^{-1}Q_{h}\|_\scalH \\
	     = & \sup_{u\in {\cal S}_{h}} \inf_{w \in {\cal S}_{h}} 
	      \frac{\|B^{*}A_{0}^{-1}u-w\|_\scalH}{\|u\|_\scalH} \\
	     = & \sup_{v\in {\cal S}_{h}} \inf_{w \in {\cal S}_{h}} 
	  \frac{\|B^{*}v-w\|_\scalH}{\|A_{0}v\|_\scalH} 	 
\end{align*}
and as a consequence, $\overset{\circ}{\varepsilon}_\scalH(h) \leq  
{\cal O}(\overset{\circ}{\gamma}(h))$ as $h\rightarrow 0. \quad \blacksquare$

\subsection{Elliptic boundary value problems; Spectral Methods}
Let $\Omega$ be the unit square in ${\Bbb R}^{2}:\ [0,\, 1]\times [0,\, 1]$
and let there be given a function $\bb =[b_{1}(\bx),\, b_{2}(\bx)]
\in C^{1}(\Omega)\times C^{1}(\Omega)$
such that $\bb=0$ on $\partial \Omega$.  Consider the differential 
operator given by
$$
A(\bx,\, D)u=-\Delta u + \bb(\bx)\cdot\nabla u +c(\bx)u  \quad  \mbox{in $\Omega$}
$$
with $u=0$ on $\partial \Omega$.  The regularity results of Kadlec 
\cite{kadlec}, 
for example, show that $Dom(A)=H^{2}\cap H^{1}_{0}$.

The application of a Fourier Galerkin method involves trial functions 
of the form $\phi_{\bk}(\bx)=\sin(k_{1}\pi x_{1})\cdot \sin(k_{2}\pi 
x_{2})$.
If  $|\bk|=k_{1}+k_{2}$ denotes the length of the 
 multi-index $\bk$,  parameterize the family of subspaces as
 ${\cal S}_{h}=\mbox{span} 
\begin{Sb}
	|\bk|h<1
\end{Sb}
 \{\phi_{\bk}\}$, and assign ${\cal S}_{h}={\cal S}_{1,h}={\cal S}_{2,h}$.  
 Define $A_{0}=-\Delta$ with $Dom(A_{0})=H^{2}(\Omega)\cap 
 H^{1}_{0}(\Omega)$ and $Dom(A_{0}^{1/2})={\cal V}=H^{1}_{0}(\Omega)$.
 Observe that $A_{0}$ is positive definite and self-adjoint in ${\cal 
 H}=L^{2}(\Omega)$ and
$$
 A_{0}\phi_{\bk}=\lambda^{0}_{\bk} \phi_{\bk},
$$
with $\lambda^{0}_{\bk}=(k_{1}^{2}+k_{2}^{2})\pi^{2}$.
If $B$ denotes the closure of $\bb(\bx)\cdot\nabla \, + c(\bx)$ on 
$C_{0}^{\infty}(\Omega)$ then the ${\cal H}$-adjoint of $B$ may be 
calculated as
$B^{*}u=- \nabla\cdot(\bb(\bx) u)=-\bb(\bx)\cdot\nabla u+(c(\bx)-\nabla\cdot\bb(\bx))u$
and $Dom(B^{*})\supset Dom(A_{0})$.
Since $B^{*}$ is compact relative to $A_{0}$ , $B^{*}$ has relative 
$A_{0}$-bound of $0$ and Theorem \ref{Hcase} Part \ref{RelCmpt} asserts that 
$\overset{\circ}{\varepsilon}_\scalH(h)\rightarrow 0$.  Furthermore, Theorem 
\ref{Hcase} Part \ref{RelBnded} 
provides a mechanism for estimating the rate at which 
$\overset{\circ}{\varepsilon}_\scalH(h)\rightarrow 0$.

First, notice that if $\{\lambda,\, {\hat u} \}$ is an eigenpair for $A$, then
$$
-\Delta {\hat u}= (\lambda - c){\hat u} - \bb\cdot\nabla {\hat u} \in 
H_{0}^{1}(\Omega)=Dom(A_{0}^{1/2}).
$$
Thus, ${\hat u}\in Dom(A_{0}^{3/2})$ and
\begin{align*}
	\|(I-Q_{h}){\hat u}\|^{2}= & \sum_{|\bk|h\geq 1} |\langle 
	           \phi_{\bk},\ {\hat u}\rangle|^{2} \\
	           = & \sum_{|\bk|h\geq 1} |\langle 
	            A_{0}^{-3/2}\phi_{\bk},\ A_{0}^{3/2}{\hat u}\rangle|^{2} \\
	          = & \sum_{|\bk|h\geq 1} (\lambda^{0}_{\bk})^{-3}|\langle 
	            \phi_{\bk},\ A_{0}^{3/2}{\hat u}\rangle|^{2} \\
	          \leq  & h^{6}\left(\frac{2}{\pi^{2}}\right)^{3}\sum_{|\bk|h\geq 1} |\langle 
	            \phi_{\bk},\ A_{0}^{3/2}{\hat u}\rangle|^{2} \\
	            \leq  &  h^{6}\left(\frac{2}{\pi^{2}}\right)^{3}
	            \| A_{0}^{3/2}{\hat u}\|_\scalH^{2}.
\end{align*}
The first inequality is a consequence of $k_{1}^{2}+k_{2}^{2}\geq 
\frac{1}{2}|\bk|^{2}\geq 1/(2h^{2})$; the second is Bessel's inequality 
with respect to the orthonormal system $\{\phi_{\bk}\}$.

A similar argument can be organized to estimate 
$\overset{\circ}{\gamma}(h)$.  For any $v\in L^{2}(\Omega)$, 
$A_{0}^{-1}v \in  H^{2}(\Omega) \cap H^{1}_{0}(\Omega) $ and
$B^{*}A_{0}^{-1}v \in  H^{1}_{0}(\Omega)  \subset Dom(A_{0}^{1/2})$.
  $A_{0}^{1/2}B^{*}A_{0}^{-1}$ is a closed, everywhere 
defined operator on $ L^{2}(\Omega)$, hence $A_{0}^{1/2}B^{*}A_{0}^{-1}$ is a bounded operator
on $ L^{2}(\Omega)$.

 Now for any $v\in {\cal H}$ with 
$\|v\|_\scalH=1$, 
\begin{align*}
	\inf_{w \in {\cal S}_{h}} \|B^{*}A_{0}^{-1}v-w\|_\scalH^{2}
	      = & \sum_{|\bk|h\geq 1} |\langle 
		           \phi_{\bk},\ B^{*}A_{0}^{-1}v\rangle|^{2} \\
		           = & \sum_{|\bk|h\geq 1} |\langle 
		            A_{0}^{-1/2}\phi_{\bk},\ A_{0}^{1/2}B^{*}A_{0}^{-1}v\rangle|^{2} \\
		          = & \sum_{|\bk|h\geq 1} (\lambda^{0}_{\bk})^{-1}|\langle 
		            \phi_{\bk},\ A_{0}^{1/2}B^{*}A_{0}^{-1}v\rangle|^{2} \\
		          \leq  & h^{2}\left(\frac{2}{\pi^{2}}\right)
		          \sum_{|\bk|h\geq 1} |\langle 
		            \phi_{\bk},\ A_{0}^{1/2}B^{*}A_{0}^{-1}v\rangle|^{2} \\
		            \leq  & h^{2}\left(\frac{2}{\pi^{2}}\right)
		            \| A_{0}^{s}B^{*}A_{0}^{-1}v\|_\scalH^{2} \leq 
		            M h^{2}
\end{align*}  
for some $M <\infty$.
From this we obtain
$$
	\overset{\circ}{\gamma}(h) \leq  \sup_{v\in {\cal H}} \inf_{w \in {\cal S}_{h}} 
		  \frac{\|B^{*}A_{0}^{-1}v-w\|_\scalH}{\|v\|_\scalH} \leq {\cal O}(h).
$$

Theorem \ref{MainH_unb} asserts that there exists an $h_{0}>0$ sufficiently small 
so that for each $h<h_{0}$ and all $u\in {\cal U}$, there is a 
$u_{h} \in {\cal U}_{h}$ so that
$$
	\|u_{h}-Q_{h}u\|_\scalH \leq c\, h^{4},
	$$
whereas $\|u_{h}-u\|_\scalH$ and $\|u-Q_{h}u\|_\scalH$ will each 
be only of order $h^{3}$ in general.

\section{Application to Krylov Subspace Methods}

  Our previous results and examples were formulated in an 
  infinite-dimensional setting and had (ultimately) an asymptotic 
  character as the dimensions of the approximating subspaces grew 
  without bound.  Although, strictly speaking, such results are not applicable 
  in finite-dimensional settings (e.g., if $A$ is an $n \times n$ 
  matrix),   asymptotic results of this sort {\em can} be useful in understanding
the characteristics of methods for calculating eigenvalues of very
large matrices   that are expected to give accurate results
with low dimension approximating subspaces.  We proceed with an analysis of
the biorthogonal Lanczos method and of the Arnoldi method
with this understanding in mind.  Throughout
this section ${\cal H}={\Bbb C}^n$, equipped with the usual Euclidean norm.

  Let $\bA$ be an $n \times n$ complex
	matrix and let $\bJ$ be a tridiagonal matrix that is 
	similar to $\bA$ -- so that $\bA=\bV \bJ \bV^{-1}$ for some $n \times n$ 
	invertible matrix $\bV$.  For any index $1 \leq \ell \leq n $, let 
	$\bJ_{\ell}$ denote the $\ell$th principal submatrix of $\bJ$:
	$$
     \bJ_{\ell}=\left[	\begin{array}{cccccc}
		\alpha_{1} & \gamma_{2} &  &  &  &   \\
		\beta_{2} & \alpha_{2} & \gamma_{3} &  &  &   \\
		 & \beta_{3} & \alpha_{3} & \ddots &  &   \\
		 &  & \ddots & \ddots &  &   \\
		 &  &  &  &  & \gamma_{\ell}  \\
		 &  &  &  & \beta_{\ell} & \alpha_{\ell}
	\end{array}  \right]
	$$
 and define $\tilde{\bJ}_{\ell+1}$ via a partitioning of $\bJ$ as
	$$
	\bJ=\left[ 
	\begin{array}{cc}
		\bJ_{\ell} & \gamma_{\ell+1}\be_{\ell}\be_{1}^{t}  \\
		\beta_{\ell+1}\be_{1}\be_{\ell}^{t}  & \tilde{\bJ}_{\ell+1}
	\end{array}\right]
	$$
	Let $\bV_{\ell}=\left[ 
		\bv_{1},\, \bv_{2},\,  \ldots ,\, \bv_{\ell} \right]$ 
denote a matrix containing the first $\ell$ columns of $\bV$;
		and for $\bW^{*}=\bV^{-1}$, let $\bW_{\ell}=\left[ 
		\bw_{1},\, \bw_{2},\, \ldots ,\, \bw_{\ell} \right]$ 
denote a matrix containing the first $\ell$ columns of $\bW$.
		Denote the remaining columns of $\bV$ 
	and $\bW$ as  $\tilde{\bV}_{\ell+1}=\left[
		\bv_{\ell+1},\, \bv_{\ell+ 2},\, \ldots ,\, \bv_{n} \right]$  
		and $\tilde{\bW}_{\ell+1}=\left[ 
		\bw_{\ell+1},\, \bw_{\ell+ 2},\, \ldots ,\, \bw_{n} \right]$ 
		respectively so that,
	$\bV=\left[ 
	\begin{array}{cc}
		\bV_{\ell} & \tilde{\bV}_{\ell+1}
	\end{array}\right]$
	$\bW=\left[ 
	\begin{array}{cc}
		\bW_{\ell} & \tilde{\bW}_{\ell+1}
	\end{array}\right]$
	
	The Lanczos algorithm \cite{LanOrig} builds up the 
	matrices $\bJ$, $\bV$, and $\bW$ one column at a time starting with 
	the vectors $\bv_{1}$ and $\bw_{1}$. Only information on the 
	action of $\bA$ and $\bA^{*}$ on selected vectors in ${\Bbb C}^{n}$ 
	is used.  	Different choices for $\bv_{1}$ 
	and $\bw_{1}$ produce distinct outcomes for $\bJ$, if all goes 
	well.  Recovering from situations where not all goes well is a 
	fundamental aspect of later refinements of the algorithm 
	(though we will not dwell on it) and two approaches to this are 
	discussed in \cite{vanDoor} and \cite{PTL}.
	
	At the $\ell$th step, the basic recursion appears as
	\begin{align}
		\bA \bV_{\ell}= & 
		\bV_{\ell}\bJ_{\ell}+\beta_{\ell+1}\bv_{\ell+1}\be_{\ell}^{t} 
		\label{recur1} \\
		\bA^{*} \bW_{\ell}= & \bW_{\ell}\bJ_{\ell}^{*}+
		\gamma_{\ell+1}\bw_{\ell+1}\be_{\ell}^{t}. \nonumber
	\end{align}
Typically, normalization is determined so that $|\gamma_{i}|=|\beta_{i}|$ 
and $\bv_{i}^{*}\bw_{i}=1$.  With exact arithmetic, the first $\ell-1$ 
steps yield matrices $\bV_{\ell}$,  $\bW_{\ell}$ that satisfy 
$\bV_{\ell}^{*}\bW_{\ell}=\bI$, 
\begin{align*}
	Ran(\bV_{\ell})= & \mbox{span}\{\bv_{1},\, \bA\bv_{1},\, \dots,\, 
	\bA^{\ell-1}\bv_{1}\}={\cal K}_{\ell}(\bA,\bv_{1}) \\
	Ran(\bW_{\ell})= & \mbox{span}\{\bw_{1},\, \bA^{*}\bw_{1},\, \dots,\, 
	(\bA^{*})^{\ell-1}\bw_{1}\}={\cal K}_{\ell}(\bA^{*},\bw_{1})
\end{align*}

Using the natural parameterization provided by $\ell$ (instead of $h$), note that 
$Q_{\ell}=\bV_{\ell}\bW_{\ell}^{*}$  is a projection onto the subspace
${\cal S}_{\ell}=\mbox{span} 
\begin{Sb}
	k=1,\ldots,\ell
\end{Sb}
 \{\bv_{k}\}$ and so $A_{\ell}=Q_{\ell}\bA Q_{\ell}=\bV_{\ell}\bJ_{\ell}\bW_{\ell}^{*}$,
 and the eigenvalues for $\bJ_{\ell}$, the so-called ``Lanczos eigenvalues",
 are Galerkin approximations to the eigenvalues of $A$.  In particular, if
$\bJ_{\ell} \bz^{(\ell)}= \lambda^{(\ell)}\bz^{(\ell)} $ and
$\bu^{(\ell)}=\bV_{\ell}\bz^{(\ell)}$ then
$\{\lambda^{({\ell})},\,\bu^{(\ell)}\}$ will be a Galerkin
eigenvalue/eigenvector pair ostensibly approximating some
eigenvalue/eigenvector pair $\{\lambda,\,\bu\}$ of $\bA$.

\begin{theorem} \label{lan_thm}
	Let $\lambda$ be a simple eigenvalue of $\bA$ with an associated unit
eigenvector $\bu$. Suppose  the Lanczos eigenvalue $\lambda^{({\ell})}$ 
converges to  $\lambda$ as $\ell$ increases.
There exist constants
$c_0,\, c_{1}>0$ so that for each $\ell>1$ the
eigenvalue/eigenvector pair $\{\lambda^{({\ell})},\,\bu^{(\ell)}\}$ 
(with $\bu^{(\ell)}$ appropriately normalized)
satisfies
$$
c_0 \|\bu^{({\ell})}-Q_{{\ell}}\bu\| \leq
|\gamma_{\ell+1}|\, \|\bv_{\ell}\|\, |\bw_{\ell+1}^{*} \bu| 
\leq c_{1}\|\bu^{({\ell})}-Q_{{\ell}}\bu\| 
+|\lambda - \lambda^{({\ell})}|\, \|\bu^{({\ell})}\|
$$
Furthermore, there exists a $c>0$ such that 
	$$
	\varepsilon_{{\cal H}}(\ell)\leq 
	c\, \frac{|\bw_{\ell+1}^{*} \bu|}{|\beta_{2}\beta_{3}\cdots \beta_{\ell}|}
	$$
\end{theorem}

The proof is deferred to the end of the section.
The point here is that $\bu^{({\ell})}$ will approach $Q_{{\ell}}\bu$
at a rate governed (for the most part) by how quickly
$|\bw_{\ell+1}^{*} \bu|\rightarrow 0$.
If $\varepsilon_{{\cal H}}(\ell)\rightarrow 0$
then $\bu^{({\ell})}$ will approach $Q_{{\ell}}\bu$ more rapidly than 
it will converge to $\bu$ itself. 

    Analogous results may be obtained for the Arnoldi method.
	Suppose $\bA$ has a Schur decomposition   
	$\bA=\bV \bH \bV^{*}$ for some $n \times n$ 
	unitary matrix $\bV$ and some upper Hessenberg matrix, $\bH$.  
	For any index $1 \leq \ell \leq n $, let 
	$\bH_{\ell}$ denote the $\ell$th principal submatrix of $\bH$:
	$$
     \bH_{\ell}=\left[	\begin{array}{cccccc}
		\alpha_{1} & \gamma_{12} & \gamma_{13}  & \dots &  & \gamma_{1\ell}  \\
		\beta_{2} & \alpha_{2} & \gamma_{23} &\dots  &  &  \gamma_{2\ell} \\
		 & \beta_{3} & \alpha_{3} & \dots &  &  \gamma_{3\ell} \\
		 &  & \ddots & \ddots &  & \vdots   \\
		 &  &  &  &  & \gamma_{\ell-1,\ell}  \\
		 &  &  &  & \beta_{\ell} & \alpha_{\ell}
	\end{array}  \right]
	$$
 and define $\bH_{12}^{(\ell)}$ and $\tilde{\bH}_{\ell+1}$ 
 via a partitioning of $\bH$ as
	$$
	\bH=\left[ 
	\begin{array}{cc}
		\bH_{\ell} & \bH_{12}^{(\ell)}  \\
		\beta_{\ell+1}\be_{1}\be_{\ell}^{t}  & \tilde{\bH}_{\ell+1}
	\end{array}\right]
	$$
Partition $\bV=\left[ 
	\begin{array}{cc}
		\bV_{\ell} & \tilde{\bV}_{\ell+1}
	\end{array}\right]$ as before.  The Arnoldi method \cite{ArnOrig} 
	proceeds using a recurrence similar to (\ref{recur1}):
	$$
	\bA \bV_{\ell}=  
		\bV_{\ell}\bH_{\ell}+\beta_{\ell+1}\bv_{\ell+1}\be_{\ell}^{t} 
	$$
but here $\bV_{\ell}$ is forced to have orthonormal columns.
	
	Note that	
$Q_{\ell}=\Pi_{\ell}=\bV_{\ell}\bV_{\ell}^{*}$  is an orthogonal projection onto the 
subspace ${\cal S}_{\ell}=\mbox{span} 
\begin{Sb}
	k=1,\ldots,\ell
\end{Sb}
 \{\bv_{k}\}$.  So $A_{\ell}=Q_{\ell}\bA Q_{\ell}=\bV_{\ell}\bH_{\ell}\bV_{\ell}^{*}$,
 and the eigenvalues for $\bH_{\ell}$, the so-called ``Arnoldi eigenvalues",
 are Galerkin approximations to the eigenvalues of $\bA$.  In particular, if
$\bH_{\ell}\, \bz^{(\ell)}= \lambda^{(\ell)}\bz^{(\ell)} $ and
$\bu^{(\ell)}=\bV_{\ell}\bz^{(\ell)}$ then
$\{\lambda^{({\ell})},\,\bu^{(\ell)}\}$ will be a Galerkin
eigenvalue/eigenvector pair ostensibly approximating some
eigenvalue/eigenvector pair $\{\lambda,\,\bu\}$ of $\bA$.

The expression for $\varepsilon_{{\cal H}}(\ell)$ becomes
	$$
	\varepsilon_{{\cal H}}(\ell) =  
	\sup_{\bu\in {\cal U}}\frac{ \|\bH_{12}^{(\ell)} 
	(\tilde{\bV}_{\ell+1}^{*}\bu ) \| }
	{\|\tilde{\bV}_{\ell+1}^{*}\bu\|}
	$$	
	\begin{theorem}   \label{arn_thm}
	Let $\lambda$ be a simple eigenvalue of $\bA$ with an associated
eigenvector $\bu$. Suppose  the Arnoldi eigenvalue $\lambda^{({\ell})}$ 
converges to  $\lambda$ as $\ell$ increases.
There exist constants
$c_0,\, c_{1}>0$ so that for each $\ell>1$ the
eigenvalue/eigenvector pair $\{\lambda^{({\ell})},\,\bu^{(\ell)}\}$ 
(with $\bu^{(\ell)}$ appropriately normalized)
satisfies
$$
c_0 \|\bu^{({\ell})}-Q_{{\ell}}\bu\| \leq
  \|\bH_{12}^{(\ell)} \tilde{\bV}_{\ell+1}^{*}\bu \| 
\leq c_{1}\|\bu^{({\ell})}-Q_{{\ell}}\bu\| 
+|\lambda - \lambda^{({\ell})}|\, \|\bu^{({\ell})}\|
$$
Furthermore there exists a
$c_{2}>0$  so that
$$
	c_0 \frac{\|u^{(\ell)}-\Pi_{\ell}u\|}
	{\|(I-\Pi_{\ell})\bu\|}\leq 
	\frac{\|\bH_{12}^{(\ell)} (\tilde{\bV}_{\ell+1}^{*}\bu) \| }
	{\|\tilde{\bV}_{\ell+1}^{*}\bu\|} \leq c_{1} 
	\frac{\|u^{(\ell)}-\Pi_{\ell}\bu\|}
	{\|(I-\Pi_{\ell})\bu\|}+c_{2} 
	\|(I-\Pi_{\ell})\bu^{*}\|
$$
where $\bu^{*}$ is a left unit eigenvector associated with $\lambda$.
\end{theorem}
{\sc Proof:} An easy calculation provides
$ Q_{\ell}\bA(I-Q_{\ell})\bu=\bV_{\ell}\bV_{\ell}^{*}\bA
{\tilde \bV}_{\ell+1}{\tilde \bV}_{\ell+1}^{*}\bu $.
Both assertions then follow from (\ref{nec_eqn1}) and (\ref{nec_eqn2}) of Theorem
\ref{Necessary}. $  \quad \blacksquare $

Two preliminary results will first be established that may be of some
independent interest.
 The following two lemmas generalize Lemma 4.1 and Lemma 6.1 of 
 \cite{WatkElsn}.
 
 \begin{lemma}  \label{WatElsLemma}
	 Let ${\cal S}$ be an $s$-dimensional subspace of ${\Bbb C}^{n}$ and 
	 suppose $\bS \in {\Bbb C}^{n \times s}$ is a given matrix with 
	 orthonormal columns so that ${\cal S}=Ran(\bS)$.  Let ${\cal T}$ be 
	 a $t$-dimensional subspace of ${\Bbb C}^{n}$ with $t\geq s$.  There 
	 exists  $\bT \in {\Bbb C}^{n \times s}$ with 
	 orthonormal columns so that $Ran(\bT) \subset {\cal T}$ and
	 $$
	 \| \bS-\bT \| \leq \sqrt{2}\, \delta_{\cal H}({\cal S},{\cal T})
	 $$
 \end{lemma}
 {\sc Proof:} If ${\cal S}\cap {\cal T}$ is a nontrivial 
 subspace of ${\Bbb C}^{n}$ with dimension $k$ then there is 
 a matrix $\bS_{0}\in {\Bbb C}^{n \times k}$ with orthonormal columns 
 and a unitary matrix $\bZ_{0}\in {\Bbb C}^{s \times s}$ so that ${\cal 
 S}\cap {\cal T} = Ran(\tilde{\bS}_{0})$ and 
 $\bS\bZ_{0}=[\tilde{\bS}_{0}\, \tilde{\bS}_{1}]$ for some choice of 
 $\tilde{\bS}_{1} \in {\Bbb C}^{n \times (s-k)}$ having orthonormal 
 columns. (In particular,  if $t+s > n$ then $k \geq t+s -n$.) 
 If ${\cal S}\cap {\cal T}= \{0\}$ then define directly, $k=0$, 
 $\bZ_{0}=\bI$, and $\tilde{\bS}_{1}=\bS$.
 
 Let $\hat{\bT}_{1} \in {\Bbb C}^{n \times (t-k)}$ be a matrix with 
 orthonormal columns so that $Ran(\hat{\bT}_{1})={\cal T}\ominus ({\cal 
 S}\cap {\cal T})$ and let $\hat{\bT}_{2} \in {\Bbb C}^{n \times (n-t)}$
 be chosen so that $[\tilde{\bS}_{0}\, \hat{\bT}_{1}\, \hat{\bT}_{2}]$
 is unitary.  
 
 The orthogonal projection onto ${\cal S}$ is given by $\Pi_{\cal S}=[\tilde{\bS}_{0}\, \tilde{\bS}_{1}]
 \left[
 \begin{array}{c}
 	\tilde{\bS}_{0}^{*}  \\
 	\tilde{\bS}_{1}^{*}
 \end{array}\right]$.  The orthogonal projections onto ${\cal T}$ 
 and ${\cal T}^{\perp}$ are given by
 $$
 \Pi_{\cal T}=[\tilde{\bS}_{0}\, \hat{\bT}_{1}]
 \left[
 \begin{array}{c}
 	\tilde{\bS}_{0}^{*}  \\
 	\hat{\bT}_{1}^{*}
 \end{array}\right]  \qquad\mbox{and}\qquad
\bI-\Pi_{\cal T}=\hat{\bT}_{2}\hat{\bT}_{2}^{*},
 $$
 respectively.  Then, 
 \begin{align*}
 	\delta_{\cal H}({\cal S},{\cal T})= & \|(\bI-\Pi_{\cal T})
	 \Pi_{\cal S}\| \\
	   = & \|\hat{\bT}_{2}
	 \left(\hat{\bT}_{2}^{*}[\tilde{\bS}_{0}\, \tilde{\bS}_{1}]\right)
	 \left[
	 \begin{array}{c}
	 	\tilde{\bS}_{0}^{*}  \\
	 	\tilde{\bS}_{1}^{*}
	 \end{array}\right]   \| = \| \hat{\bT}_{2}^{*} \tilde{\bS}_{1}  \|
 \end{align*}
 
 Observe that $\left[
	 \begin{array}{c}
	 	\hat{\bT}_{1}^{*}\tilde{\bS}_{1}  \\
	 	\hat{\bT}_{2}^{*}\tilde{\bS}_{1}
	 \end{array}\right] \in {\Bbb C}^{(n-k) \times (s-k)}$ has 
	 orthonormal columns.  By simultaneously diagonalizing the summands 
	 in the expression 
	 $(\hat{\bT}_{1}^{*}\tilde{\bS}_{1})^{*}(\hat{\bT}_{1}^{*}\tilde{\bS}_{1}) 	 
	 +(\hat{\bT}_{2}^{*}\tilde{\bS}_{1})^{*}(\hat{\bT}_{2}^{*}\tilde{\bS}_{1})=\bI$,
	 one finds a unitary matrix $\bZ_{1}\in {\Bbb C}^{(s-k) \times 
	 (s-k)}$ and diagonal matrices $\Gamma=diag(\gamma_{1},\, 
	 \dots,\, \gamma_{s-k})\in {\Bbb C}^{(t-k) \times 
	 (s-k)}$, $\Sigma=diag(\sigma_{1},\, 
	 \dots,\, \sigma_{s-k})\in {\Bbb C}^{(n-t) \times 
	 (s-k)}$ with $\gamma_{1} \geq 
	 \gamma_{2}\geq \dots  \geq \gamma_{s-k}$, $\sigma_{1}\leq 
	 \sigma_{2}\leq\dots \leq \sigma_{s-k}$, and 
	 $\sigma_{i}^{2}+\gamma_{i}^{2}=1$ so that
	 \begin{align*}
	 	\bZ_{1}^{*}(\hat{\bT}_{1}^{*}\tilde{\bS}_{1})^{*}(\hat{\bT}_{1}^{*}
		 \tilde{\bS}_{1})\bZ_{1} & = \Gamma^{*}\Gamma \\
		 \bZ_{1}^{*}(\hat{\bT}_{2}^{*}\tilde{\bS}_{1})^{*}
		 (\hat{\bT}_{2}^{*}\tilde{\bS}_{1})\bZ_{1} & = \Sigma^{*}\Sigma.
	 \end{align*}
Thus there are unitary matrices $\bU_{1}\in {\Bbb C}^{(t-k) 
\times(t-k)}$ and $\bU_{2}\in {\Bbb C}^{(n-t) \times (n-t)}$ such that
$\hat{\bT}_{1}^{*}\tilde{\bS}_{1}\bZ_{1}=\bU_{1}\Gamma$ and 
$\hat{\bT}_{2}^{*}\tilde{\bS}_{1}\bZ_{1}=\bU_{2}\Sigma$.  Thus, 
$\delta_{{\cal H}}({\cal S},{\cal T})=
\| \hat{\bT}_{2}^{*} \tilde{\bS}_{1}  \| =\| \Sigma  \| 
=\sigma_{s-k}$.

Now let $\bU_{11}$ denote the $(t-k) \times (s-k) $ matrix consisting 
of the first $s-k$ columns of $\bU_{1}$ and define 
$$
\bT=\left[ \tilde{\bS}_{0} \ 
\hat{\bT}_{1}\bU_{11}\bZ_{1}^{*}\right]\, \bZ_{0}^{*}
$$
Then 
\begin{align*}
	\|\bS-\bT\|^{2} = & \|(\bS-\bT)\bZ_{0}\, diag(\bI,\bZ_{1}) \|^{2}\\
	          = & \| [\tilde{\bS}_{0}\, \tilde{\bS}_{1}\bZ_{1}] - 
	          [\tilde{\bS}_{0}\, \hat{\bT}_{1}\bU_{11}] \|^{2} 
	           =   \| \tilde{\bS}_{1}\bZ_{1}- \hat{\bT}_{1}\bU_{11}\|^{2} \\
	      =& \max_{\|\bx\|=1} 
	      \bx^{*}(\tilde{\bS}_{1}\bZ_{1})^{*}(\tilde{\bS}_{1}\bZ_{1}) \bx+
	     \bx^{*}(\hat{\bT}_{1}\bU_{11})^{*}(\hat{\bT}_{1}\bU_{11}) \bx \\
            & \qquad - 2 \Re e\,\bx^{*}
             \bU_{11}^{*}(\hat{\bT}_{1}^{*}\tilde{\bS}_{1}\bZ_{1})\bx \\
	     = & 2 \max_{\|\bx\|=1} (1-\Re e\,\sum_{i=1}^{s-k}\gamma_{i}\, 
	     |x_{i}|^{2}) = 2(1-\gamma_{s-k})
\end{align*}
Since $1+\gamma_{s-k}>1$, the conclusion follows from
$$
\|\bS-\bT\|\leq 
\sqrt{2(1-\gamma_{s-k})(1+\gamma_{s-k})}=\sqrt{2}\,\sigma_{s-k}\ \quad \blacksquare
$$

\begin{lemma} \label{LanczosGap}
	Let ${\cal U}$ be the invariant subspace of $\bA$ associated with the 
	eigenvalue $\lambda$.  Suppose $\ell$ is an index such that $\ell 
	\geq m=\dim\,{\cal U}$.  Then
	$$
	|\beta_{2}\beta_{3}\cdots \beta_{\ell+1}|\leq 
	\|\bW\|(1+\sqrt{2})\delta_{{\cal H}}({\cal U},{\cal 
	K}_{\ell}(\bA,\bv_{1}))
	$$
\end{lemma}
{\sc Proof:} Let $\Pi_{\ell}$ denote the orthogonal projection onto ${\cal 
	K}_{\ell}(\bA,\bv_{1})$.  Fix for the time being a vector $\bu\in {\cal 
	U}$ with $\|\bu\|=1$ and let $p$ be that polynomial of degree 
	$\ell-1$ or less such that $p(A)\bv_{1}=\Pi_{\ell}u$.  Pick 
	$\epsilon >0$ and (if necessary) vary the coefficients of $p$ 
	slightly to obtain a polynomial $\hat{p}$ so that 
	$\hat{p}(\lambda)\neq 0$ and 
	$\|(\bI-\Pi_{\ell})\bu\|+\epsilon\geq \|\bu - 
	\hat{p}(\bA)\bv_{1}\|$.  Now, $\hat{p}(\bA)$ maps ${\cal U}$ 
	bijectively to itself so there exists $\hat{\bu}\in {\cal U}$ such 
	that $\hat{p}(\bA)\hat{\bu}=\bu$.  Then
$$
\|(\bI-\Pi_{\ell})\bu\|+\epsilon\geq \|\bu - 
	\hat{p}(\bA)\bv_{1}\|=\|\hat{p}(\bA)(\hat{\bu} -\bv_{1})\|.
$$

Denote the $QR$ decomposition of $ \tilde{\bW}_{\ell+1}$ as $ 
\tilde{\bW}_{\ell+1}=\bS\bR$ so that $\bS$ has orthonormal columns 
and $Ran( \tilde{\bW}_{\ell+1})=Ran(\bS)$.  Observe that $Ran( 
\tilde{\bW}_{\ell+1})=(Ran(\bV_{\ell}))^{\perp}={\cal 
	K}_{\ell}(\bA,\bv_{1})^{\perp}$.  Lemma \ref{WatElsLemma} guarantees 
	the existance of a matrix $\bT\in {\Bbb C}^{n \times (n-\ell)}$ with 
	orthonormal columns and with $Ran(\bT)\subset {\cal U}^{\perp}$ such 
	that
$$
	 \| \bS-\bT \| \leq \sqrt{2}\, \delta_{{\cal H}}({\cal 
	K}_{\ell}(\bA,\bv_{1})^{\perp},{\cal U}^{\perp})=\sqrt{2}\, \delta_{{\cal 
	H}}({\cal U} ,{\cal K}_{\ell}(\bA,\bv_{1})
 $$
Note that
\begin{align*}
	\tilde{\bW}_{\ell+1}^{*}\hat{p}(\bA)\bv_{1}= & \tilde{\bW}_{\ell+1}^{*}
	(\bV\hat{p}(\bA)\bW^{*})\bV\be_{1} \\
	       = & [0\, \bI_{n-\ell}]\hat{p}(\bJ)\be_{1} \\
	       = & (\beta_{2}\beta_{3}\cdots \beta_{\ell+1})\be_{1}\in {\Bbb C}^{(n-\ell)}
\end{align*}

On the other hand, 
$$
\tilde{\bW}_{\ell+1}^{*}\hat{p}(\bA)\bv_{1} = \bR^{*}(\bS-\bT)^{*}\hat{p}(\bA)\bv_{1}
     +\bR^{*}\bT^{*}\hat{p}(\bA)(\hat{\bu}-\bv_{1}) -\bR^{*}\bT^{*}\hat{p}(\bA)\hat{\bu}
$$
But since $\bT^{*}\hat{p}(\bA)\hat{\bu}=0$, we find
\begin{align*}
	|\beta_{2}\beta_{3}\cdots \beta_{\ell+1}| = & 
	\|\tilde{\bW}_{\ell+1}^{*}\hat{p}(\bA)\bv_{1}\| \\
	 \leq & \|\bR\| (\|\bS-\bT\|\| \hat{p}(\bA)\bv_{1} \|+ 
	 \|\hat{p}(\bA)(\hat{\bu}-\bv_{1})\|\\
	 \leq & \|\tilde{\bW}_{\ell+1}\| (\sqrt{2}\,\delta_{{\cal H}}({\cal U} ,{\cal K}_{\ell}(\bA,\bv_{1})) + \|(\bI-\Pi_{\ell})\bu\|+\epsilon) 
\end{align*}
The conclusion follows by letting $\epsilon \rightarrow 0$, taking 
the supremum over all unit vectors $\bu \in {\cal U}$, and noticing 
that $\|\tilde{\bW}_{\ell+1}\|\leq \|\bW\|. \quad \blacksquare $

{\sc Proof of Theorem \ref{lan_thm}:} An easy calculation provides
$ Q_{\ell}\bA(I-Q_{\ell})\bu=\gamma_{\ell+1}\bv_{\ell}\bw_{\ell+1}^{*} \bu $.
The first assertion then follows from (\ref{nec_eqn1}) of Theorem
\ref{Necessary}.  For the second assertion, observe first that
since $\lambda$ is simple
 $\|(\bI-\Pi_{\ell})\bu\|=\delta_{{\cal H}}({\cal U} ,
{\cal K}_{\ell}(\bA,\bv_{1})$.  Then using Lemma \ref{LanczosGap} we 
find,
\begin{align*}
	\varepsilon_{{\cal H}}(\ell)= & \frac{\|Q_{\ell}\bA(I-Q_{\ell})\bu\|}
	{\|(\bI-\Pi_{\ell})\bu\|}\\
   \leq  & \|\bW\|(1+\sqrt{2})
\frac{|\gamma_{\ell+1}|\,\|\bv_{\ell}\|\,|\bw_{\ell+1}^{*} \bu|}
	                  {|\beta_{2}\beta_{3}\cdots \beta_{\ell+1}|}\\
	               & \leq c\,\frac{|\bw_{\ell+1}^{*} 
	               \bu|}{|\beta_{2}\beta_{3}\cdots \beta_{\ell}|},
\end{align*}
since $|\gamma_{\ell+1}|=|\beta_{\ell+1}|.  \quad \blacksquare $

\section*{Appendix: Lower Bounds to {\em sep}}  \label{prelim}
\renewcommand{\thesection}{\Alph{section}}
      \setcounter{section}{1}
      \setcounter{theorem}{0}
      \setcounter{equation}{0}
 Let ${\cal W}_{1}$ and ${\cal W}_{2}$ be complex Hilbert spaces and 
 denote with ${\cal L}({\cal W}_{2},{\cal W}_{1})$ the associated Banach space 
 of bounded linear transformations from ${\cal W}_{2}$ to ${\cal W}_{1}$.
Suppose there are given two linear 
operators, $L_{1}:{\cal W}_{1}\rightarrow {\cal W}_{1}$, $L_{2}:{\cal 
W}_{2}\rightarrow {\cal W}_{2}$ such that $L_{1}$ is closed, densely 
defined (but not necessarily bounded) on ${\cal W}_{1}$ 
and $L_{2}$ is bounded (and everywhere defined) in ${\cal W}_{2}$.
Define an operator $\mathcal{T}:{\cal L}({\cal 
W}_{2},Dom(L_{1}))\rightarrow {\cal L}({\cal W}_{2},{\cal W}_{1})$ as
$\mathcal{T}(S)=L_{1}S-SL_{2}$ and let
\begin{equation}
	sep(L_{1},L_{2})\stackrel{\rm def}{=} \inf_{S\in {\cal L}({\cal 
	W}_{2},Dom(L_{1}))} \frac{\|\mathcal{T}(S)\|_{{\cal W}_{2}\rightarrow {\cal W}_{1}}}
	{\|S\|_{{\cal W}_{2}\rightarrow {\cal W}_{1}}},
	\label{sepdef}
\end{equation}
so that, in particular, if $S$ solves $\mathcal{T}(S)=M$ and 
$0<\eta<sep(L_{1},L_{2})$ then 
$$
\|S\|_{{\cal W}_{2}\rightarrow {\cal W}_{1}} \leq 
	\frac{1}{sep(L_{1},L_{2})}
	\|M\|_{{\cal W}_{2}\rightarrow {\cal W}_{1}}\leq 
	\frac{1}{\eta}
	\|M\|_{{\cal W}_{2}\rightarrow {\cal W}_{1}}.
$$
The following results 
are mild generalizations of \cite{Rosenblum} Theorem 3.1, p. 264 
and \cite{Heinz} Theorem 5, p. 427.
For all $\epsilon>0$, we define the pseudospectral sets
\begin{align*}
	\Lambda_{\epsilon}(L)= & \left\{ z \in {\Bbb C}\left| \ 
	\|(z-L)^{-1}\|\ge \frac{1}{\epsilon} \right. \right\}.
\end{align*}
As $\epsilon \rightarrow 0$, $\Lambda_{\epsilon}(L)$ shrinks to  
$\sigma(L)$.

\begin{lemma} \label{sepbnd}
 Suppose that $\sigma(L_{1})$ and $\sigma(L_{2})$ are disjoint.
Then for every $M \in {\cal L}({\cal W}_{2},{\cal W}_{1})$, the operator equation
\begin{equation}	
	L_{1}S-SL_{2}=M
	\label{sylvester}
\end{equation}
has a unique solution $S$ that is a bounded linear transformation 
from  ${\cal W}_{2}$ to $Dom(L_{1})\subset{\cal W}_{1}$.  That is, 
$\mathcal{T}(S)$ is bijective.

(a) $S$ has the representation
\begin{equation}
	S=\frac{1}{2\pi \imath}
	\int_{\Gamma_{2}}(L_{1}-z)^{-1}M(z-L_{2})^{-1}\ dz.
	\label{Srep}
\end{equation}
where $\Gamma_{2}$ consists of a finite number of closed 
rectifiable Jordan curves enclosing all points of $\sigma(L_{2})$ and 
no points of $\sigma(L_{1})$. 

(b) Let $\epsilon_{1},\ 
\epsilon_{2}>0$ be chosen so that 
$\Lambda_{\epsilon_{1}}(L_{1})\cap {\Gamma}_{2} = \emptyset$ 
and $\Lambda_{\epsilon_{2}}(L_{2})\cap
{\Gamma}_{2} = \emptyset$, respectively. That is, 
$\Lambda_{\epsilon_{1}}(L_{1})$ lies entirely {\em outside} 
${\Gamma}_{2}$ and $\Lambda_{\epsilon_{2}}(L_{2})$ lies entirely 
{\em inside} ${\Gamma}_{2}$.
Then 
\begin{equation}
	sep(L_{1},L_{2}) \geq  
	\frac{2 \pi\epsilon_{1}\epsilon_{2}}{\text{length}(\Gamma_{2})}
	\label{sylbound1}
\end{equation}
where $\text{length}(\Gamma_{2})$ is the arc length of $\Gamma_{2}$.

(c) If the numerical ranges of $L_{1}$ and $L_{2}$, denoted
respectively ${\frak w}(L_{1})$ and ${\frak w}(L_{2})$, 
are disjoint sets in ${\Bbb C}$ 
and  $\Delta=dist({\frak w}(L_{1}),\  {\frak w}(L_{2})) > 0$, then
\begin{equation}
		sep(L_{1},L_{2}) \geq   \Delta
	\label{sylbound2}
\end{equation}
\end{lemma}
{\sc Proof:} A trivial extension of Corollary 3.3 of \cite{Rosenblum} guarantees that 
 there are a finite number of closed 
rectifiable Jordan curves enclosing all points of $\sigma(L_{2})$ and 
no points of $\sigma(L_{1})$.  We orient each of the curves 
positively (i.e., counterclockwise) and refer to them collectively as $\Gamma_{2}$.
Observe that ${\Gamma}_{2}\subset \rho(L_{1})\cap 
\rho(L_{2})$, so the right hand side of \ref{Srep} is 
well-defined and maps ${\cal W}_{2}$ into $Dom(L_{1})\subset{\cal W}_{1}$. 
Since $\frac{1}{2\pi \imath}
		\int_{\Gamma_{2}}(z-L_{2})^{-1}\ dz \ = I$, $\frac{1}{2\pi \imath}
		\int_{\Gamma_{2}}(z-L_{1})^{-1}\ dz \ = 0$, and
		$L(L-z)^{-1}=I+z(L-z)^{-1}$, direct substitution of \ref{Srep} 
into \ref{sylvester} yields
\begin{align*}
	L_{1}S-SL_{2}&=\frac{1}{2\pi \imath}
		\int_{\Gamma_{2}}M(z-L_{2})^{-1}-(L_{1}-z)^{-1}M\ dz \\
		& =M \left(\frac{1}{2\pi\imath}\int_{\Gamma_{2}}(z-L_{2})^{-1}dz\right)
		-\left(\frac{1}{2\pi\imath}\int_{\Gamma_{2}}(L_{1}-z)^{-1}\ dz\right)M \\
		& = M.
\end{align*}
Suppose $N$ is the difference between any two solutions of 
\ref{sylvester}.  Then $L_{1}N=NL_{2}$ and
\begin{align*}
	(z-L_{1})N & =N(z-L_{2}) \\
	N(z-L_{2})^{-1} & =(z-L_{1})^{-1}N \\
	N\left(\frac{1}{2\pi 
			\imath}\int_{\Gamma_{2}}(z-L_{2})^{-1}dz\right) & =
			\left(\frac{1}{2\pi \imath}\int_{\Gamma_{2}}(z-L_{1})^{-1}dz\right)N \\
	  N & = 0.
\end{align*}
Thus \ref{Srep} gives the one, unique solution to \ref{sylvester}.
  
To show (b) note that,
  \begin{align*}
	  \|S\|_{{\cal W}_{2}\rightarrow {\cal W}_{1}} & \leq \frac{1}{2\pi}
	\int_{\Gamma_{2}}\|(L_{1}-z)^{-1}\|_{{\cal W}_{1}}
	\|M\|_{{\cal W}_{2}\rightarrow {\cal W}_{1}}\|(z-L_{2})^{-1}\|_{{\cal 
	W}_{2}}\ |dz|     \\
	& \leq \frac{\text{length}(\Gamma_{2})}{2\pi}\max_{z\in \Gamma_{2}}\|(L_{1}-z)^{-1}\|_{{\cal W}_{1}}
	\|M\|_{{\cal W}_{2}\rightarrow {\cal W}_{1}}
	\max_{z\in \Gamma_{2}}\|(z-L_{2})^{-1}\|_{{\cal W}_{2}}	  \\
	& \leq \frac{\text{length}(\Gamma_{2})}{2\pi}\frac{1}{\epsilon_{1}}
	\cdot 	\frac{1}{\epsilon_{2}}\|M\|_{{\cal W}_{2}\rightarrow {\cal W}_{1}}  
  \end{align*}

For (c), note that existence and uniqueness of $S$ follows from (a) and the 
observation that disjoint numerical ranges of $L_{1}$ and $L_{2}$ imply 
disjoint spectra for $L_{1}$ and $L_{2}$, however, we will use a 
different representation for $S$ to obtain estimates.
Define ${\hat z}_{1},\ {\hat z}_{2}\in {\Bbb C}$ so that 
${\hat z}_{1} \in closure\{ {\frak w}(L_{1})\} $, $ {\hat z}_{2}\in 
closure\{ {\frak w}(L_{2})\} $, and
\begin{align}
\Delta = |{\hat z}_{1}- {\hat z}_{2}| = &  \inf \begin{Sb}
 z_{1}\in {\frak w}(L_{1}) \\
 z_{2}\in {\frak w}(L_{2})
 \end{Sb}  |z_{1}-z_{2}| \\
 = & \inf \begin{Sb}
 x \in {\cal W}_1 \\
 y \in {\cal W}_2
 \end{Sb}  |\frac{\langle x,\ L_{1}x\rangle_{{\cal W}_1}}{\langle x,\ x\rangle_{{\cal W}_1}}
 - \frac{\langle y,\ L_{2}y\rangle_{{\cal W}_2}}{\langle y,\ y\rangle_{{\cal W}_2}}|.
\end{align}
Define $\theta \in [0,\ 2\pi)$ as $\theta = arg({\hat z}_{1}- {\hat z}_{2} ) $. 
Let ${\hat L_{1}}=e^{-\imath \theta}(L_{1} - {\hat z}_{1})$ and 
${\hat L_{2}}=e^{-\imath \theta}(L_{2} - {\hat z}_{1})$.
Then ${\frak w}({\hat L_{1}})=e^{-\imath \theta}({\frak w}(L_{1}) - {\hat z}_{1})$ and 
${\frak w}({\hat L_{2}})=e^{-\imath \theta}({\frak w}(L_{2}) - {\hat z}_{1})$.

The goal of this translation and rotation is to insure
$$
\Re e\ \langle x,\ {\hat L_{1}}x\rangle_{{\cal W}_1}
 \geq 0
$$
for all $x\in \ Dom(L_{1})$, and 
$$
\Re e \ \langle y,\ {\hat L_{2}}y\rangle_{{\cal W}_2}
\leq -\Delta 
\|y\|_{{\cal W}_2}^{2}
$$
for all $y\in \ {\cal W}_2)$.

Under these circumstances, the Lumer-Phillips theorem (cf. \cite{pazy}, 
\S1.4) guarantees that $-{\hat L}_{1}$ and ${\hat L}_{2}$ each generate strongly continuous, 
one parameter semigroups, $e^{-t{\hat L}_{1}}$ and $e^{t{\hat L}_{2}}$, respectively. Furthermore, 
$\|e^{-t{\hat L}_{1}}\|_{{\cal W}_1}\leq 1$ and 
$\|e^{t{\hat L}_{2}}\|_{{\cal W}_2}\leq e^{-t\Delta}$ for all $t\geq 0$.

Notice that $S$ satisfies (\ref{sylvester}) if and only if it is also a 
solution to $	{\hat L_{1}}S-S{\hat L_{2}}=e^{-\imath \theta}M $ and 
this leads to the following representation for $S$:
\begin{equation}
	S=e^{-\imath \theta}\int_{0}^{\infty}e^{-t{\hat L}_{1}}Me^{t{\hat L}_{2}}\ dt.
	\label{Srep2}
\end{equation}
Indeed, note that with this expression for $S$, we have for any $v\in {\cal W}_{2}$
	\begin{align*}
	\|\left(L_{1}S-S L_{2}-M\right)v &\|_{{\cal W}_{1}} = 
	\|\left({\hat L_{1}}S-S{\hat L_{2}}-e^{-\imath \theta}M\right) v\|_{{\cal W}_{1}} \\
		= & \|\int_{0}^{\infty}{\hat L_{1}}e^{-t{\hat L}_{1}}Me^{t{\hat L}_{2}}v
		-e^{-t{\hat L}_{1}}Me^{t{\hat L}_{2}}{\hat L_{2}}v\ dt -Mv \|_{{\cal W}_{1}}\\
 = & \|\int_{0}^{\infty}-\frac{d}{dt}\left(e^{-t{\hat L}_{1}}Me^{t{\hat L}_{2}}v\right)\ dt -Mv\|_{{\cal W}_{1}}\\
 =	& \|	\lim_{t\rightarrow 0}\left(e^{-t{\hat L}_{1}}Mv-Mv\right)+
		\lim_{t\rightarrow 0}e^{-t{\hat L}_{1}}M\left(e^{t{\hat L}_{2}v-v}v\right) \\
		& \quad -\lim_{t\rightarrow \infty}
		\left(e^{-t{\hat L}_{1}}Me^{t{\hat L}_{2}}v\right) \|_{{\cal W}_{1}} \\
 \leq & \lim_{t\rightarrow 0}\|(e^{-t{\hat L}_{1}}Mv-Mv\|_{{\cal W}_{1}}+
		\|M\|_{{\cal W}_{2}\rightarrow {\cal W}_{1}}
		\lim_{t\rightarrow 0}\|e^{t{\hat L}_{2}}v-v\|_{{\cal W}_{1}} \\
		& \quad +\lim_{t\rightarrow \infty}
		\|M\|_{{\cal W}_{2}\rightarrow {\cal W}_{1}}e^{-t\Delta}\|v\|_{{\cal 
		W}_{2}} = 0.	
		\end{align*}
Immediately then,  one obtains
\begin{align}
	\|S\|_{{\cal W}_{2}\rightarrow {\cal W}_{1}}  & \leq 
	\int_{0}^{\infty}\|e^{-t{\hat L}_{1}}\|_{{\cal W}_{1}} \|M\|_{{\cal W}_{2}\rightarrow {\cal 
	W}_{1}}\|e^{t{\hat L}_{2}}\|_{{\cal W}_{2}}\ dt \\
	& \leq \int_{0}^{\infty} e^{-t\Delta}\ dt \ 
	\|M\|_{{\cal W}_{2}\rightarrow {\cal W}_{1}}
	 = \left(\frac{1}{\Delta}\right)\|M\|_{{\cal W}_{2}\rightarrow {\cal W}_{1}} 
	 	\label{bound2}
\end{align}


\begin{thebibliography}{9}
 
 	\bibitem{ArnOrig}  Arnoldi, W. E. (1951), The principle of 
 	minimized iterations in the solution of the matrix eigenvalue problem,  
 	{\it Quart. Appl. Math.} \textbf{9}, pp. 17--29.  
    
   	\bibitem{BabOsb}  Babu\v{s}ka, I., and J. Osborn (1991), Eigenvalue 
 	Problems, in ``\textit{Finite Element Methods}'' Handbook of Numerical Analysis, Vol. 2,
    edited by P. G. Ciarlet and J. L. Lions.  Elsevier Science 
    Publisher (North Holland)
    
 	\bibitem{BabOsb2}  Babu\v{s}ka, I., and J. Osborn (1989), Finite 
 	element--Galerkin approximation of the eigenvalues and eigenvectors 
 	of selfadjoint problems,  {\it Math. Comp.} \textbf{52}, pp. 275--297.  
 	
 	\bibitem{Chat}  Chatelin, F. (1983), \textit{Spectral Approximation of 
 	Linear Operators}, (Academic Press, New York)
 	
 	\bibitem{Des1978}  Descloux, J., N. Nassif, and J. Rappaz (1978), 
 	On spectral approximation, Part I: The problem of convergence, {\it 
 	RAIRO Anal. Num\'{e}r.} \textbf{12}, pp. 97--112. 
 	
 	\bibitem{Des1981}  Descloux, J., M. Luskin, and J. Rappaz (1981), 
 	Approximation of the spectrum of closed operators: the determination 
 	of normal modes of a rotating basin, {\it 
 	Math. Comp.} \textbf{36}, pp. 137--154. 
 	
 \bibitem{vanDoor} Grimme, E. J., D. C. Sorensen, and P. van Dooren, (1995)
         Model reduction of state space systems via an implicitly restarted
         Lanczos method.
 	 
  \bibitem{Heinz} Heinz, E. (1951), Beitr\"{a}ge zur St\"{o}rungstheorie der 
 	 Spektralzerlegung, {\it Math. Ann.} \textbf{123}, pp. 415--438.
 	 
 	 \bibitem{kadlec} Kadlec, J. (1964), The regularity of the solution
of the Poisson problem in a domain whose boundary is similar to that of a
convex domain, {\it Czech Math J. } \textbf{14}(89), pp. 386--393. (In Russian)
 	 
 	 \bibitem{Kato0} Kato, T. (1960), Estimation of iterated matrices, 
 	with application to the von Neumann condition, {\it Numerische 
 	Mathematik} \textbf{2}, pp. 22--29.

    \bibitem{Kato} Kato, T. (1976), \textit{ Perturbation Theory for Linear 
    Operators}, (Springer, Heidelberg)	
 
    \bibitem{Kny}  Knyazev, A. (1994),  New estimates for Ritz vectors, CIMS 
 	NYU Tech Report 677 (to appear in Mathematics of Computation)
 	
   \bibitem{LanOrig}  Lanczos, C. (1950), An iteration method for the 
 	solution of the eigenvalue problem of linear differential and 
 	integral operators,  {\it J. Res. Nat. Bur. Stand.} \textbf{45}, pp. 255--282.  
 	
    \bibitem{PTL} Parlett, B. N., D. R. Taylor, and Z. A. Liu, A 
    look-ahead Lanczos algorithm for unsymmetric matrices, {\it Math. Comp.} 
    \textbf{44}, pp. 105--124.
    
    \bibitem{pazy} Pazy, A. (1983), \textit{Semigroups of 
 	Linear Operators and Applications to Partial Differential Equations},
 	 (Springer, Heidelberg)	
    
    \bibitem{Rosenblum} Rosenblum, M. (1956), On the operator equation 
    $BX-XA=Q$, {\it Duke Math. J.} pp. 263--269.
    
    \bibitem{WatkElsn} Watkins, D. S., and L. Elsner (1991), Convergence
of algorithms of decomposition type for the eigenvalue problem, 
     {\it Lin. Alg. Appl.} \textbf{143} pp. 19--47.

    \bibitem{wloka} Wloka, J. (1987), \textit{Partial Differential 
    Equations}, (Cambridge Univ. Press, Cambridge)
    
 \end{thebibliography}
\end{document}